\documentclass[draft,preprint,12pt,numbers,sort&compress]{elsarticle}
\usepackage{mathrsfs}
\usepackage{amssymb}
\usepackage{amsthm}
\usepackage{mathrsfs}
\usepackage[centertags]{amsmath}
\usepackage{amsfonts}

\usepackage{color}

\input amssym.def
\input amssym.tex

\date{}

\newtheorem{Theorem}{Theorem}[section]

\newtheorem{Lemma}{Lemma}[section]

\newcommand\R{\mbox{\bf R}}

\newcommand\SR{\mbox{\scriptsize\bf R}}

\newcommand{\definition}{{\lower .5ex
  \hbox{$\>\>\stackrel{\triangle}{=}\>\>$} }}
\newcommand\supp{\mathop{\rm supp}}

\begin{document}

\baselineskip=22pt
\thispagestyle{empty}
\baselineskip=22pt
\thispagestyle{empty}

\mbox{}
\bigskip

\begin{center}{\Large\bf The Cauchy problem for the Ostrovsky equation with positive dispersion}\\[1ex]
{Wei Yan\footnote{Email: yanwei19821115@sina.cn}$^{a,d}$,
\quad Yongsheng Li \footnote{Email: yshli@scut.edu.cn}$^b$,
\quad Jianhua  Huang\footnote{ Email: jhhuang32@nudt.edu.cn}$^c$,
\quad Jinqiao  Duan\footnote{ Email: duan@iit.edu}$^{d}$}\\[1ex]

{$^a$College of Mathematics and Information Science, Henan Normal University,}\\
{Xinxiang, Henan 453007,   China}\\[1ex]

{$^b$School of Mathematics, South  China  University of  Technology,}\\
{Guangzhou, Guangdong 510640,  China}\\[1ex]

{$^c$College of Science, National University of Defense  Technology,}\\
{ Changsha, Hunan 410073,  China}\\[1ex]

{$^d$ Department of Applied Mathematics, Illinois Institute of Technology,}\\[1ex]
{Chicago, IL 60616, USA  }\\[1ex]

\end{center}

\noindent{\bf Abstract.}This paper is devoted to studying the Cauchy problem
for the Ostrovsky equation
\begin{eqnarray*}
      \partial_{x}\left(u_{t}-\beta \partial_{x}^{3}u
  +\frac{1}{2}\partial_{x}(u^{2})\right)
    -\gamma u=0,
\end{eqnarray*}
with positive $\beta$ and $\gamma $. This equation describes the propagation
 of surface waves in a rotating  oceanic flow. We first prove that the
  problem is  locally well-posed
in $H^{-\frac{3}{4}}(\R)$. Then  we  reestablish the bilinear estimate, by means of  the
  Strichartz estimates  instead of calculus inequalities and  Cauchy-Schwartz
   inequalities.  As a byproduct,  this bilinear estimate leads to the proof of the
    local well-posedness of the problem in $H^{s}(\R)$ for $ s>-\frac{3}{4}$, with help of a fixed point argument.

\bigskip
\bigskip

{\large\bf 1. Introduction}
\bigskip

\setcounter{Theorem}{0} \setcounter{Lemma}{0}

\setcounter{section}{1}

We consider the Cauchy problem for the Ostrovsky equation  with positive dispersion
\begin{eqnarray}
&&\partial_{x}\left(u_{t}-\beta \partial_{x}^{3}u
  +\frac{1}{2}\partial_{x}(u^{2})\right)
    -\gamma u=0,    \;\; \mbox{ with } \beta>0, \gamma>0,            \label{1.01}\\
    &&u(x,0)=u_{0}(x),  \label{1.02}
\end{eqnarray}
where $u$ represents the free surface of a liquid and the positive parameter $\gamma$
measures the effect of rotation.
This equation (\ref{1.01})   was proposed by Ostrovsky  \cite{O} as a model
for weakly nonlinear long waves in a rotating liquid, by taking into account of the Coriolis force.
It describes the propagation of surface waves in the ocean in a rotating frame of reference.
In fact,   $\beta$ determines the type of dispersion, more precisely, $\beta<0$ (negative dispersion)
 for  surface and internal waves in the ocean
or surface waves in a shallow channel with an uneven bottom,  and $\beta>0$ (positive
dispersion) for capillary waves on the surface of a liquid or for oblique magneto-acoustic
waves in plasma \cite{Be,GaSt,GiGrSt,Gri,L1981,S1981,S1986}).

The Ostrovsky equation  (\ref{1.01})  can be rewritten in the following form
\begin{equation}
   u_{t}-\beta\partial_{x}^{3}u
    +\frac{1}{2}\partial_{x}(u^{2})- \gamma \partial_{x}^{-1} u
    =0.\label{1.03}
\end{equation}
When $\gamma = 0$,
it  reduces to  the Korteweg-de Vries equation which has been investigated widely
 \cite{B,CKSTT, Bourgain-GAFA93,Bou,CCT,KPV,KPV1991,KPV1993,KPV19931,KPV0,T,Tzvetkov}.
 By introducing the Fourier restriction norm method, Bourgain \cite{Bourgain-GAFA93} proved that
 the Cauchy problem for the KdV equation is globally well-posed in $L^{2}$ for the periodic case and
 nonperiodic case.  By using the  Fourier  restriction   norm method,
 Kenig et al. \cite{KPV,KPV0} proved that the Cauchy problem for the KdV equation
  is locally well-posed in $H^{s}(\R),s>-\frac{3}{4}$ and ill-posed in $H^{s}(\R)$
   with $s<-\frac{3}{4}$ in the   sense that the data-to-solution map fails
   to be uniformly continuous as map from $H^{s}(\R)$ to $C_{t}^{0}H^{s}(\R)$.
   This means that $s=-\frac{3}{4}$ is
    the critical  regularity index  in Sobolev spaces for KdV equation.
By using the I-method and the Fourier restriction norm method, Colliander et al. \cite{CKSTT}
 showed  that the Cauchy
    problem for the KdV equation is globally well-posed in $H^{s}(\R)$ with $s>-\frac{3}{4}$.
Guo \cite{G} and Kishimoto \cite{Kis}  proved that   the Cauchy  problem for the KdV equation  is
globally well-posed in  $H^{-\frac{3}{4}}(\R)$ with  the aid of the $I$-method and
the modified Besov spaces.  Kappeler and  Topalov \cite{KT} proved that the flow map
 extends continuously
to $H^{-1}$ in the periodic case with the aid of  inverse scattering transformation.
 Molinet \cite{Molinet}
 proved that no well-posedness result can possibly hold below $s=-1$ in the periodic
  case
  in the sense that the solution map
 of KdV equation does not extend to a continuous map  from $H^{s}$ for $s<-1$ to
  distribution.
 Molinet \cite{Molinet2011} proved that the solution-map  associated with the KdV
equation cannot be continuously extended in $H^{s}(R)$ for $s<-1$.
Liu \cite{L}  established a priori bounds
for KdV equation in $H^{-\frac{3}{4}}(\R)$. Buckmaster and  Koch \cite{BK} proved the
existence of weak solutions to the KdV initial value problem on the real line with $H^{-1}$
 initial data and  studied  the problem of orbital and asymptotic $H^{s}$ stability of
 solitons for integers $s=-1;$
 and  established new a priori $H^{-1}$ bound for solutions to the KdV equation.

The stability of the solitary waves or soliton solutions of  the Ostrovsky equation
(\ref{1.01}) has also been examined
 \cite{LLSIAM,LLDCDS,LL2007,LV,L2007,VL2005,ZL}.
Choudhury et al. \cite{CIL} studied the
Hamiltonian formulation, nonintegrability and local bifurcations for the Ostrovsky equation.
Others have studied the Cauchy problem for (\ref{1.01}); for instance,
see \cite{GL, GH,HJ,HJ0, I,LLDCDS, IM1,IM2,IM3,IM4,IM5,LM,LV,Molinet,VL,WC,Z,ZL}.
The results in  \cite{IM1,IM3,Tsu} showed that $s=-\frac{3}{4}$
is the critical regularity index for (\ref{1.01}) in Sobolev spaces.   Recently,
Coclite and Ruvo \cite{CR2014,CR2015}
 have investigated the convergence of the Ostrovsky equation to the Ostrovsky-Hunter equation,
  and also the dispersive and diffusive limits for the Ostrovsky-Hunter type equation.
Moreover, Li et al. \cite{LHY} proved that the Cauchy problem for the Ostrovsky
equation with negative dispersion is locally well-posed in $H^{-\frac{3}{4}}(\R).$
However, the well-posedness  of the Ostrovsky equation with positive dispersion  in $H^{-\frac{3}{4}}(\R)$
has not yet been   shown up to now.
\medskip

Observe that if $u(x,t)$ is the solution to the Cauchy problem for (\ref{1.03}),
then $u^{\lambda}(x,t)=\lambda^{-2}u\left(\frac{x}{\lambda},\frac{t}{\lambda^{3}}\right)$,
for $\lambda > 0$,
is the solution to the following    equation
\begin{eqnarray}
&&u^{\lambda}_{t}-\beta\partial_{x}^{3}u^{\lambda}+\frac{1}{2}\partial_{x}((u^{\lambda})^{2})-\gamma\lambda^{-4}
\partial _{x}^{-1}u^{\lambda}=0\label{1.04},\\
&&u^{\lambda}(x,0)=\lambda^{-2}u_{0}\left(\frac{x}{\lambda}\right)\label{1.05}.
\end{eqnarray}
If $u$ is the solution to (\ref{1.03}), then $v(x,t)=\beta^{-1}u(x,\beta^{-1}t)$ is the solution $v_{t}-v_{xxx}+\frac{1}{2}\partial_{x}(v^{2})-\beta^{-1}\gamma\partial_{x}^{-1}v=0$. Hence without loss of generality,
we can assume that $\gamma=\beta=1$ in this paper.

\medskip

In the present paper,  we first prove that
 (\ref{1.03}) with initial condition (\ref{1.02})  is locally well-posed in $H^{-\frac{3}{4}}(\R)$.  More precisely, we establish a
bilinear estimate with $s=-\frac{3}{4}$,  which combines  with the fixed point theorem to yield the local well-posedness of the
 Cauchy problem for (\ref{1.03}) in $H^{-\frac{3}{4}}(\R)$.
Then by using the Strichartz estimates instead
of     calculus inequalities and  the
     Cauchy-Schwartz inequality, we  reestablish the bilinear estimate  for the Ostrovsky
      equation (\ref{1.03}) with $s>-\frac{3}{4},$  which
       combines with Lemma 2.8 (below) and the fixed point theorem
      to imply the local well-posedness of  the
 Cauchy problem for (\ref{1.03}) in $H^{s}(\R)$ with $s>-\frac{3}{4}.$

\bigskip

Before stating the main results, we introduce some notations.
 Throughout this paper, we assume that $\lambda\geq1$ and
$C$ is a positive constant
 which may vary from line to line. Let  $\epsilon$ be a small number with
 $0<\epsilon<{10^{-4}}$. Note that
$a\sim b$ means that $|b|\leq |a|\leq 4|b|$,
 $a\gg b$ means that $|a|> 4|b|.$  Let $\psi $ be  a smooth function with
 support in $[-1,2]$ and taking the value
 $1$ in $[-1,1]$.  For a set $I\subset \R^{2}$,  define the indicator function $\chi_{I}((\xi,\tau))=1$ if $(\xi,\tau)\in I$;
  $\chi_{I}((\xi,\tau))=0$ if $(\xi,\tau)  \notin I$. Let   $\mathscr{F}u$ be the Fourier
  transformation of $u$ with respect to
 both space and time variables,  and $\mathscr{F}^{-1}u$ be the corresponding  inverse
 transformation, while $\mathscr{F}_{x}u$ denotes
  the Fourier transformation of $u$
 with respect to the space variable and $\mathscr{F}^{-1}_{x}u$ denotes the corresponding
 inverse transformation.
Define
\begin{eqnarray*}
&&\langle \cdot\rangle=1+|\cdot|,\\
&&D^{\prime}:=\left\{(\xi,\tau)\in \R^{2}:|\xi|\leq\frac{1}{8}, |\tau|\geq |\xi|^{-3}\right\},\\
&&D_{1}:=\left\{(\xi,\tau)\in \R^{2}:|\xi|\leq\frac{1}{8}, |\tau|< |\xi|^{-3}\right\},\\
&&D_{2}:=\left\{(\xi,\tau)\in \R^{2}:\frac{1}{8}<|\xi|\leq1, |\tau|< |\xi|^{-3}\right\},\\
&&D_{3}:=\left\{(\xi,\tau)\in \R^{2}:\frac{1}{8}<|\xi|\leq1, |\tau|\geq |\xi|^{-3}\right\},\\
&&A_{j}:=\left\{(\xi,\tau)\in \R^{2}: 2^{j}\leq \langle\xi\rangle<2^{j+1}\right\},\\
&&\phi^{\lambda}(\xi)=\xi^{3}+\frac{1}{\lambda^{4}\xi},\sigma^{\lambda}=\tau+\phi^{\lambda}(\xi),
\sigma_{j}^{\lambda}=\tau_{j}+\phi^{\lambda}(\xi_{j})(j=1,2),\\
&&B_{k}:=\left\{(\xi,\tau)\in \R^{2}: 2^{k}\leq \left\langle\sigma^{\lambda}
\right\rangle<2^{k+1}\right\},\\
&&S^{\lambda}(t)\phi=e^{t(\partial_{x}^{3}+\lambda^{-4}\partial_{x}^{-1})}\phi
=C\int_{\SR}e^{-it (\xi^{3}+\xi^{-1}\lambda^{-4})}\mathscr{F}_{x}\phi(\xi)d\xi.
\end{eqnarray*}
 Obviously, $\left\{(\xi,\tau)\in \R^{2}:|\xi|\leq1, \tau\in \R\right\}=D^{\prime}\cup D_{1}\cup D_{2}\cup D_{3}.$
 Here $j,k$ are nonnegative integers.
 Space
$
  X_{\lambda}^{s, \>b}
$ is defined by
$$
X_{\lambda}^{s, \>b}= \left\{u\in  \mathscr{S}^{'}(\R^{2}) :\, \|u\|_{X_{\lambda}^{s, \>b}}
 =  \left\|\langle\xi\rangle^{s} \left\langle\sigma^{\lambda}\right\rangle^{b}\mathscr{F}u(\xi,\tau)
 \right\|_{L_{\tau\xi}^{2}(\SR^{2})}<\infty\right\}.
$$
$X_{\lambda}^{s,b}$  was
 introduced by  Rauch and Reed \cite{RR}, Beals \cite{Beals},
 Bourgain \cite{Bourgain-GAFA93},
 Klainerman and Machedon \cite{KM},  and further developed by Kenig, Ponce and Vega \cite{KPV1993}.
Space
$
  X_{\lambda}^{s,\,b,\> 1} = \left\{u\in  \mathscr{S}^{'}(\R^{2})\, :\,
   \|u\|_{X_{\lambda}^{s,\>b,\>1}}<\infty\right\}
$ is equipped with the following  norm
\begin{eqnarray*}
    \|u\|_{X_{\lambda}^{s,\>b,\>1}}
&=& \biggl\|\biggl(\left\|\langle\xi\rangle^{s}
    \left\langle\sigma^{\lambda}\right\rangle ^{b} \mathscr{F}u(\xi,\tau)
     \right\|_{L_{\tau\xi}^{2}(A_{j}\cap B_{k})}\biggr)_{j,\>k\geq 0}
      \biggr\|_{\ell_{j}^{2}(\ell_{k}^{1})}
       \nonumber\\
&\sim& \biggl[\sum_{j}2^{2js}\biggl(\sum_{k}2^{bk}\|\mathscr{F}u(\xi,\tau)\|_{L_{\tau\xi}^{2}
        (A_{j}\cap B_{k})}\biggr)^{2}\biggr]^{1/2}.
\end{eqnarray*}
Space $ X_{\lambda} $ is defined by
\begin{eqnarray*}
&&X_{\lambda}\nonumber\\&& = \left\{u\in  \mathscr{S}^{'}(\R^{2}) :\, \|u\|_{X_{\lambda}}=\left\|\mathscr{F}^{-1}
 [\chi_{D^{c}}\mathscr{F}u]\right\|_{X_{\lambda}^{-\frac{3}{4}, \frac{1}{2}, 1}}
    +\|\mathscr{F}^{-1}[\chi_{ D^{\prime}}\mathscr{F}u]\|_{X_{\lambda}^{-\frac{3}{4}, \frac{1}{2}}}<\infty\right\}
\end{eqnarray*}
and space $Y$ is defined by
 $
 Y= \left\{u\in  \mathscr{S}^{'}(\R^{2})\, :\, \|u\|_{Y}=
 \left\|\langle\xi\rangle^{-\frac34}\mathscr{F}u(\xi,\tau)\right\|_{L_{\xi}^{2}L_{\tau}^{1}}<\infty\right\},
$
where $D^{c}\cup D^{\prime}=\R^2_{\tau\xi}$, respectively. Spaces
 $\hat{X_{\lambda}}, \hat{X_{\lambda}}^{s,\>b,\>1}, \hat{X_{\lambda}}^{s,\>b}$ are equipped with the following norms
\begin{eqnarray*}
\|f\|_{\hat{X_{\lambda}}}=\|\mathscr{F}^{-1}f\|_{X_{\lambda}},\quad
\|f\|_{\hat{X_{\lambda}}^{s,b,1}}= \|\mathscr{F}^{-1}f\|_{X_{\lambda}^{s,b,1}},
\quad
\|f\|_{\hat{X_{\lambda}}^{s,b}}=\|\mathscr{F}^{-1}f\|_{X_{\lambda}^{s,b}},
\end{eqnarray*}
respectively.
The space $ X_{_{\lambda},T}$ denotes the restriction of $X_{\lambda}$ onto the finite time interval $[-T,T]$ and
is equipped with the norm
 \begin{equation*}
    \|u\|_{X_{_{\lambda},\>T}} =\inf \left\{\|w\|_{X_{\lambda}}:w\in X_{\lambda}, u(t)=w(t)
 \>\> {\rm for} \> -T\leq t\leq T\right\}.
 \end{equation*}

The main results of this paper are stated in Theorem 1.1 and Theorem 1.2.

\begin{Theorem}\label{Thm1} (Well-posedness in $H^{-\frac{3}{4}}(\R)$)
The Cauchy problem  for (\ref{1.03}) is locally well-posed in $H^{-\frac{3}{4}}(\R)$. That is,
for $u_{0} \in H^{-\frac{3}{4}}(\R)$, there exists a  positive number $T$ such that the solution map $u_0\mapsto u(t)$ is locally
Lipschitz continuous from $ H^{-\frac{3}{4}}(\R)$ into $C([-T, T]; H^{-\frac{3}{4}}(\R))\cap  \|u\|_{X_{_{\lambda},T}}$.
Moreover, the solution to the Cauchy problem for (\ref{1.03}) on the time interval $[-T,T]$ is unique in the space  $X_{\lambda,T}.$
\end{Theorem}
\noindent {\bf Remark 1:}
Isaza and  J. Mej\'{\i}a \cite{IM3} and Taugawa \cite{Tsu} have proved that the Cauchy problem
  (\ref{1.03}) - (\ref{1.02})  is locally well-posed in $H^{s}(\R)$
with $s>-\frac{3}{4}$.  Moreover, Isaza and  J. Mej\'{\i}a \cite{IM3}  have proved that
 the  problem  (\ref{1.03})- (\ref{1.02})  is not  quantitatively well-posed in
   $H^{s}(\R)$ with $s<-\frac{3}{4}$. Thus, $s=-\frac{3}{4}$ is the critical
   regularity index in Sobolev space for (\ref{1.03})- (\ref{1.02}).

\noindent {\bf Remark 2:} Inspired by \cite{IK,IKT,Kis} and in view of the structure of the
Ostrovsky equation with positive dispersion,  we  choose the space  $X_{\lambda}$,
which is slightly  different from space $X$ of \cite{Kis}.
More precisely,  Kishimoto \cite{Kis} used  the set
\begin{eqnarray*}
D:=\left\{(\xi,\tau)\in \R^{2},|\xi|\leq
1,|\tau|\geq |\xi|^{-3}\right\}
\end{eqnarray*}
in \cite{Kis}.  However,
we need to take
\begin{eqnarray*}
D^{\prime}:=\left\{(\xi,\tau)\in \R^{2},|\xi|\leq \frac{1}{8},|\tau|\geq |\xi|^{-3}\right\}.
\end{eqnarray*}
 Now we give a specific example to explain the reason why we take  $D^{\prime}.$
In proving  (vii) of Lemma 3.1 (below) for the case $2^{k_{2}}=2^{k_{\rm max}}>4{\rm max}\left\{2^{k},
2^{k_{1}}\right\}$ and $(\tau_{2},\xi_{2})\in D^{\prime}$, we have   $C|\xi|^{-3}\leq 2^{k_{2}}
\sim |\sigma^{+}-
\sigma_{1}^{+}-\sigma_{2}^{+}|\sim 2^{k_{2}}=|3\xi\xi_{1}\xi_{2}-\frac{\xi_{1}^{2}
+\xi_{1}\xi_{2}+\xi_{2}^{2}}{\lambda^{4}\xi\xi_{1}\xi_{2}}|\sim C|\xi\xi_{1}\xi_{2}|
\sim |\xi_{2}|2^{2j_{1}}$, and  $|\xi_{2}|\geq C2^{-\frac{j_{1}}{2}}$,
which is crucial in establishing (vii) of Lemma 3.1.
But on the set $D=\left\{(\tau,\xi)\in \R^{2},|\xi|\leq 1,|\tau|\geq |\xi|^{-3}\right\}$,
 we can not  guarantee  that
$|\xi_{2}|\geq C2^{-\frac{j_{1}}{2}}$.

\noindent {\bf Remark 3:} Now we  outline  the proof for Theorem 1.1.
 Comparing with  the KdV equation and the Ostrovsky equation with the negative
dispersion,
the structure of the Ostrovsky equation with positive dispersion  is much  more
complicated.
More precisely, for $\lambda > 0$,  note that
\begin{equation}
  (\xi_{1}+\xi_{2})^{3}+\frac{1}{\lambda^{4} (\xi_{1}+\xi_{2})}
  - \xi_{1}^{3}-\frac{1}{\lambda ^{4}\xi_{1}}-\xi_{2}^{3}-\frac{1}{\lambda^{4}\xi_{2}}
  = 3\xi_{1}\xi_{2}(\xi_{1}+\xi_{2})-\frac{\xi_{1}^{2}+\xi_{1}\xi_{2}+\xi_{2}^{2}}
  {\lambda^{4}\xi_{1}\xi_{2}(\xi_{1}+\xi_{2})},    \label{1.06}
\end{equation}
\begin{eqnarray}
&& \xi_{1}^{3}+\frac{1}{\lambda ^{4}\xi_{1}}+ \xi_{2}^{3}+\frac{1}{\lambda ^{4}\xi_{2}}
    -\frac{(\xi_{1}+\xi_{2})^{3}}{4}-\frac{4}{\lambda ^{4}(\xi_{1}+\xi_{2})}
 \nonumber\\&&= \frac{3}{4}(\xi_{1}+\xi_{2})(\xi_{1}-\xi_{2})^{2}\left[1+
  \frac{4}{3\lambda^{4}\xi_{1}\xi_{2}(\xi_{1}+\xi_{2})^{2}}\right],   \label{1.07}
\end{eqnarray}
and
\begin{eqnarray}
&&\frac{\xi_{2}^{3}}{4}+\frac{4}{\lambda ^{4}\xi_{2}}
-\left((\xi_{1}+\xi_{2})^{3}+\frac{1}{\lambda^{4}(\xi_{1}+\xi_{2})}\right)
+\left(\xi_{1}^{3}+\frac{1}{\lambda^{4}\xi_{1}}\right)\nonumber\\&&
=\frac{3}{4}\xi_{2}(2\xi_{1}+\xi_{2})^{2}
\left[1-\frac{4}{3\lambda^{4}(\xi_{1}+\xi_{2})\xi_{1}\xi_{2}^{2}}\right].    \label{1.08}
\end{eqnarray}
From (\ref{1.06})-(\ref{1.08}),  we know that there exist no   positive
 constants $C_{j} (j=1,2,3)$ such that
\begin{eqnarray}
 &&\hspace{-1cm}\left|(\xi_{1}+\xi_{2})^{3}+\frac{1}{\lambda^{4}
 (\xi_{1}+\xi_{2})}- \xi_{1}^{3}-\frac{1}{\lambda ^{4}\xi_{1}}-\xi_{2}^{3}
 -\frac{1}{\lambda^{4}\xi_{2}}\right|\nonumber\\&&=
 \left|3\xi_{1}\xi_{2}(\xi_{1}+\xi_{2})-\frac{\xi_{1}^{2}+\xi_{1}\xi_{2}+\xi_{2}^{2}}
  {\lambda^{4}\xi_{1}\xi_{2}(\xi_{1}+\xi_{2})}\right|\nonumber\\&&
  \geq C_{1}|\xi_{1}\xi_{2}(\xi_{1}+\xi_{2})|,  \label{1.09}
\end{eqnarray}
\begin{eqnarray}
&& \left|\xi_{1}^{3}+\frac{1}{\lambda ^{4}\xi_{1}}+ \xi_{2}^{3}+\frac{1}{\lambda ^{4}\xi_{2}}
    -\frac{(\xi_{1}+\xi_{2})^{3}}{4}-\frac{4}{\lambda ^{4}(\xi_{1}+\xi_{2})}\right|
 \nonumber\\&&=\left| \frac{3}{4}(\xi_{1}+\xi_{2})(\xi_{1}-\xi_{2})^{2}\left[1+
  \frac{4}{3\lambda^{4}\xi_{1}\xi_{2}(\xi_{1}+\xi_{2})^{2}}\right]\right|
  \geq C_{2}|\xi_{1}\xi_{2}(\xi_{1}+\xi_{2})|,   \label{1.010}
\end{eqnarray}
and
\begin{eqnarray}
&&\left|\frac{\xi_{2}^{3}}{4}+\frac{4}{\lambda ^{4}\xi_{2}}
-\left((\xi_{1}+\xi_{2})^{3}+\frac{1}{\lambda^{4}(\xi_{1}+\xi_{2})}\right)
+\left(\xi_{1}^{3}+\frac{1}{\lambda^{4}\xi_{1}}\right)\right|\nonumber\\&&
=\left|\frac{3}{4}\xi_{2}(2\xi_{1}+\xi_{2})^{2}
\left[1-\frac{4}{3\lambda^{4}(\xi_{1}+\xi_{2})\xi_{1}\xi_{2}^{2}}\right]\right|
\geq C_{3}|\xi_{1}\xi_{2}(\xi_{1}+\xi_{2})|\label{1.011}
\end{eqnarray}
for any $\xi_{1},\xi_{2}\in \R$.
For the KdV  equation, we have that
\begin{eqnarray}
 &&\left|(\xi_{1}+\xi_{2})^{3}
  - \xi_{1}^{3}-\xi_{2}^{3}\right|
  = 3\left|\xi_{1}\xi_{2}(\xi_{1}+\xi_{2})\right|,\label{1.012}\\
  &&\left|\xi_{1}^{3}+ \xi_{2}^{3}
    -\frac{(\xi_{1}+\xi_{2})^{3}}{4}\right|
 = \left|\frac{3}{4}(\xi_{1}+\xi_{2})(\xi_{1}-\xi_{2})^{2}\right|,\label{1.013}\\
 &&\left|\frac{\xi_{2}^{3}}{4}
-(\xi_{1}+\xi_{2})^{3}
+\xi_{1}^{3}\right|=\left|\frac{3}{4}\xi_{2}(2\xi_{1}+\xi_{2})^{2}\right|,\label{1.014}
\end{eqnarray}
respectively. The equalities
(\ref{1.013})-(\ref{1.014}) ensure that Lemmas 3.2, 3.3 of \cite{Kis} are valid. The estimate
(\ref{1.09}) and Lemmas 3.2, 3.3 of \cite{Kis} are the key tools of  proving the main result
of Theorem 1.2 of \cite{Kis}. From Lemmas 2.3, 2.6 of \cite{LHY} and the following inequality
\begin{equation}
  (\xi_{1}+\xi_{2})^{3}-\frac{1}{\lambda^{4} (\xi_{1}+\xi_{2})}
  - \xi_{1}^{3}+\frac{1}{\lambda ^{4}\xi_{1}}-\xi_{2}^{3}+\frac{1}{\lambda^{4}\xi_{2}}
  = 3\xi_{1}\xi_{2}(\xi_{1}+\xi_{2})+\frac{\xi_{1}^{2}+\xi_{1}\xi_{2}+\xi_{2}^{2}}
  {\lambda^{4}\xi_{1}\xi_{2}(\xi_{1}+\xi_{2})},\label{1.015}
\end{equation}
we know  that the proof of Theorem 1.1 of  \cite{LHY}  is similar to the
 proof of Theorem 1.2 of \cite{Kis}.

\noindent From (\ref{1.06})-(\ref{1.011}),  we know that we cannot completely
 follow the method of \cite{Kis, LHY}
to establish the local well-posedness of the Cauchy problem for the Ostrovsky equation
with positive dispersion in $H^{-\frac34}(\R).$
Thus, for $\lambda \geq 1$ and the Ostrovsky equation with  positive  dispersion,
 we consider cases
\begin{eqnarray}
&&\left|3\xi_{1}\xi_{2}(\xi_{1}+\xi_{2})-\frac{\xi_{1}^{2}+\xi_{1}\xi_{2}+\xi_{2}^{2}}
{\lambda^{4}\xi_{1}\xi_{2}(\xi_{1}+\xi_{2})}\right|<
\frac{|\xi_{1}\xi_{2}(\xi_{1}+\xi_{2})|}{4},\label{1.016}\\
&&\left|3\xi_{1}\xi_{2}(\xi_{1}+\xi_{2})-\frac{\xi_{1}^{2}+\xi_{1}\xi_{2}+\xi_{2}^{2}}
{\lambda^{4}\xi_{1}\xi_{2}(\xi_{1}+\xi_{2})}\right|
\geq\frac{|\xi_{1}\xi_{2}(\xi_{1}+\xi_{2})|}{4},\label{1.017}
\end{eqnarray}
\begin{eqnarray}
&&\left|1+
  \frac{4}{3\lambda^{4}\xi_{1}\xi_{2}(\xi_{1}+\xi_{2})^{2}}\right|
  <\frac{1}{2},\label{1.018}\\
  &&\left|1+
  \frac{4}{3\lambda^{4}\xi_{1}\xi_{2}(\xi_{1}+\xi_{2})^{2}}\right|
  \geq\frac{1}{2},\label{1.019}
\end{eqnarray}
and
\begin{eqnarray}
&&\left|1-\frac{4}{3\lambda^{4}(\xi_{1}+\xi_{2})\xi_{1}\xi_{2}^{2}}\right|
<\frac{1}{2},\label{1.020}\\
&&\left|1-\frac{4}{3\lambda^{4}(\xi_{1}+\xi_{2})\xi_{1}\xi_{2}^{2}}\right|
\geq\frac{1}{2}.\label{1.021}
\end{eqnarray}
Only when (\ref{1.017}), (\ref{1.019}) and  (\ref{1.021}) are valid simultaneously,
 can we follow the method of \cite{Kis}.
Other cases have to be handled with the aid of techniques introduced in the present paper.

\medskip

Our second main result is stated in Theorem 1.2.

\begin{Theorem}\label{Thm2} (Bilinear estimate for  $s>-\frac{3}{4}$)\\
Let $s>-\frac{3}{4}$ and $\beta>0$ and $\gamma>0$. Then the following bilinear estimate holds:
\begin{eqnarray}
\left\|\partial_{x}(u_{1}u_{2})\right\|_{X_{\lambda}^{s,-\frac{1}{2}+2\epsilon}}\leq C\prod_{j=1}^{2}\|u_{j}\|_{X_{\lambda}^{s,\frac{1}{2}+\epsilon}}.\label{1.022}
\end{eqnarray}
\end{Theorem}

\noindent {\bf Remark 4:} By  using  the
  calculus inequalities and  Cauchy-Schwartz
   inequalities,
Isaza and Mej\'{\i}a \cite{IM1,IM3}  and Tsugawa \cite{Tsu} have already established Theorem 1.2.
But in this paper, we use the Strichartz estimates
to present an alternative proof of  Theorem 1.2.
As a byproduct,  Theorem 1.2,  combining with Lemma 2.8 (below) and a fixed point argument,
   leads to the well-posedness of  the Cauchy problem for the
Ostrovsky equation with positive dispersion
  in $H^{s}(\R)$  with $s>-\frac{3}{4}$.

\noindent {\bf Remark 5:}
Now we  outline  the proof of  Theorem 1.2.
 Isaza and  Mej\'{\i}a \cite{IM1, IM3} have used certain  calculus inequalities
 and Cauchy-Schwartz inequality
to establish this crucial bilinear estimate,  which plays the key role in proving
the local well-posedness of the Cauchy problem for the Ostrovsky equation
 with    positive dispersion in $H^{s}(\R)$ with $s>-\frac{3}{4}$.
In this paper, we use the Strichartz estimates instead of   calculus inequalities
 and Cauchy-Schwartz inequality to reestablish the same bilinear estimate
    for the Ostrovsky equation
 with   positive dispersion.
From Lemma 2.9 (below), we know that
\begin{eqnarray}
&&\left\|\int_{\tiny\begin{array}{cc}\xi=\sum\limits_{j=1}^{2}\xi_j\\ \tau
=\sum\limits_{j=1}^{2}\tau_j\end{array}}|\xi_{1}^{2}-\xi_{2}^{2}|^{s}
\prod_{j=1}^{2}\mathscr{F}u_{j}(\xi_{j},\tau_{j})
d\xi_{1}d\tau_{1}\right\|_{L_{\xi\tau}^{2}}\nonumber\\&&\leq
\left\|\int_{\tiny\begin{array}{cc}\xi=\sum\limits_{j=1}^{2}\xi_j\\ \tau
=\sum\limits_{j=1}^{2}\tau_j\end{array}}
\left|\xi_{1}^{2}-\xi_{2}^{2}\right|^{s}
\left|1+\frac{1}{3\lambda^{4}\xi_{1}^{2}\xi_{2}^{2}}\right|^{s}
\mathscr{F}u_{1}(\xi_{1},\tau_{1})
\mathscr{F}u_{2}(\xi_{2},\tau_{2}) d\xi_{1}d\tau_{1}\right\|_{L_{\xi\tau}^{2}}
\nonumber\\&&\leq C\prod_{j=1}^{2}
\|u_{j}\|_{X_{\lambda}^{0, \>\frac{3+2s}{4}\frac{1}{2}+\epsilon}},
\qquad 0<s\leq \frac{1}{2},\label{1.023}
\end{eqnarray}
which is useful in establishing the bilinear estimate.
When (\ref{1.017})  is valid, we can easily
 establish the bilinear estimate by using (\ref{1.023}).
When (\ref{1.016}) is valid, we have
\begin{eqnarray}
|\xi_{{\rm min}}|\sim \lambda^{-2}|\xi_{{\rm max}}|^{-1},\label{1.024}
\end{eqnarray}
where $|\xi_{{\rm min}}|:={\rm min}\left\{|\xi|,|\xi_{1}|,|\xi_{2}|\right\}$
and $|\xi_{{\rm max}}|:={\rm max}\left\{|\xi|,|\xi_{1}|,|\xi_{2}|\right\}$.
We will combine  (\ref{1.024})    and   (\ref{1.023}),    together with  a suitable
 splitting of region,   to establish the bilinear estimate.

\noindent {\bf Remark 6:} Levandosky and  Liu  \cite{LLSIAM} studied the stability of the generalized
Ostrovsky equation
\begin{eqnarray}
\left[u_{t}-\beta u_{xxx}+(u^{k})_{x}\right]_{x}=\gamma u,k\geq3,k\in N,\gamma>0,\beta\neq0\label{1.025}.
\end{eqnarray}
The techniques used in proving Theorem 1.2 are conducive to establish the  multilinear estimate
\begin{eqnarray}
\left\|\partial_{x}\prod_{j=1}^{k}(u_{j})\right\|_{X_{\lambda}^{s,-\frac{1}{2}+2\epsilon}}\leq C\prod_{j=1}^{k}\|u_{j}\|_{X_{\lambda}^{s,\frac{1}{2}+\epsilon}}\label{1.026}
\end{eqnarray}
for   $s=\frac{1}{4}$ with $k=3$ and $s=\frac{k-4}{2k}$ with $k\geq4$.
As a byproduct, (\ref{1.026}) in combination with Lemma 2.8
  yields a  result on the optimal regularity   for  (\ref{1.025}) in Sobolev spaces.

\medskip

The rest of this paper is arranged as follows.   After presenting some
preliminaries in the next section,   we establish bilinear estimate for
 $s=-\frac34$ in Lemmas 3.1-3.2 in Section 3.  Then, we prove local well-posedness   (Theorem 1.1)
  for $s=-\frac34$, and bilinear estimate (Theorem 1.2) for $s > -\frac34$, in Sections 4 and 5,  respectively.

\bigskip
\bigskip

 \noindent{\large\bf 2. Preliminaries }

\setcounter{equation}{0}

\setcounter{Theorem}{0}

\setcounter{Lemma}{0}

\setcounter{section}{2}

In this section, we establish Lemmas 2.1-2.9 which play an important role in establishing
  Lemma 3.1 and Theorems 1.1 and 1.2. More precisely, we establish Lemmas 2.1, 2.2, 2.3 and 2.6. These lemmas will be used to
  prove Lemma 3.1 which
  implies Lemma 3.2 immediately with the help of a suitable decomposition.
  Lemma 3.2 and Lemmas 2.4 - 2.5,  in combination with a fixed point argument, yield Theorem 1.1.
   Lemmas 2.7-2.9 are used to establish Theorem 1.2.

\begin{Lemma}\label{Lemma2.1} Assume  that $f_{k}(k=1,2)\in \mathscr{S}'(\R^{2})$ with
$\supp f_{k}\subset A_{j_{k}}(k=1,2)$,  and
\begin{eqnarray*}
K_{1}:=\inf\left\{|\xi_{1}-\xi_{2}|:\exists\> \tau_{1},\tau_{2},\>
s.t.\>(\xi_{k},\tau_{k})\in
\supp f_{k}(k=1,2)\right\}>0.
\end{eqnarray*}
If
\begin{eqnarray}
(\xi_{1},\tau_{1})\in
\supp f_{1}  ,\, (\xi_{2},\tau_{2})\in \supp f_{2},\>\xi_{1}\xi_{2}>0,\label{2.01}
\end{eqnarray}
or
\begin{eqnarray}
(\xi_{k},\tau_{k})\in
\supp f_{k}(k=1,2),\xi_{1}\xi_{2}<0, \quad
\left|1+\frac{4}{3\lambda ^{4}\xi^{2}\xi_{1}\xi_{2}}\right|
            > \frac{1}{2},\label{2.02}
\end{eqnarray}
then there exists a positive constant $C$
 such that the following inequalities hold
\begin{eqnarray}
&&\||\xi|^{\frac14}f_{1}\ast f_{2}\|_{L^{2}(\SR^{2})}\leq
C\prod_{k=1}^{2}\|f_{k}\|_{\hat{X_{\lambda}}^{0,\frac{1}{2},1}}
,\label{2.03}\\
&&\||\xi|^{\frac12}f_{1}\ast f_{2}\|_{L^{2}(\SR^{2})}\leq CK_{1}^{-\frac12}
\prod_{k=1}^{2}\|f_{k}\|_{\hat{X_{\lambda}}^{0,\frac{1}{2},1}}.
\label{2.04}
\end{eqnarray}
\end{Lemma}
\noindent {\bf Proof.}
Since $\tau=\tau_{1}+\tau_{2},\xi=\xi_{1}+\xi_{2},$
  by a direct computation, we have that
\begin{eqnarray}
 && \tau+\frac{\xi^{3}}{4}+
 \frac{4}{\lambda ^{4}\xi}-\left(\tau_{1}+\xi_{1}^{3}
  + \frac{1}{\lambda ^{4}\xi_{1}}\right)
  -\left(\tau_{2}+\xi_{2}^{3}+
  \frac{1}{\lambda ^{4}\xi_{2}}\right)\nonumber\\&&
= -\frac{3}{4}\xi(\xi_{1}-\xi_{2})^{2}
\left[1+\frac{4}{3\lambda ^{4}\xi^{2}\xi_{1}\xi_{2}}\right].
\label{2.05}
\end{eqnarray}
By using  (\ref{2.05}) and a proof similar to that for
 Lemma 2.1 of \cite{LHY},   we conclude  that
(\ref{2.03})-(\ref{2.04}) are valid.

This completes the proof of Lemma 2.1.

\begin{Lemma}\label{Lemma2.2} Assume that
 $f\in \mathscr{S}^{'}(\R^{2})$ , $g\in \mathscr{S}(\R^{2})$,
with
$\mathop{\rm supp} f\subset A_{j}$ for some $j\geq0$
 and $\,\Omega\subset \R^{2}$ with positive measure.
Let
\begin{eqnarray*}
K_{2}:=\inf\left\{|\xi_{1}+\xi|:\exists\, \tau,\tau_{1}\>
s.t.\> (\xi,\tau)\in \Omega,(\xi_{1},\tau_{1})\in A_j
   \right\}>0.
\end{eqnarray*}
If
\begin{eqnarray}
(\xi,\tau)\in \Omega,(\xi_{1},\tau_{1})\in \supp f,
\>\xi\xi_{1}<0,\label{2.06}
\end{eqnarray}
or
\begin{eqnarray}
(\xi,\tau)\in \Omega  ,\, (\xi_{1},\tau_{1})\in
\supp f,\>\xi\xi_{1}> 0,\quad
\left|1-\frac{4}{3\lambda ^{4}\xi\xi_{1}\xi_{2}^{2}}\right|
            > \frac{1}{2},\label{2.07}
\end{eqnarray}
then, for every $k\geq0$, there exists a positive constant  $C>0$
 such that the following inequalities hold
\begin{eqnarray}
&&\|f\ast g\|_{L^{2}(B_{k})}
\leq
C2^{\frac{k}{4}}\|f\|_{\hat{X_{\lambda}}^{0,\frac{1}{2},1}}
\||\xi|^{-\frac{1}{4}}g\|_{L^{2}(\SR^{2})}
\label{2.08},\\
    &&\|f\ast g\|_{L^{2}(\Omega\cap B_{k})}
\leq  C2^{\frac{k}{2}}K_{2}^{-\frac{1}{2}}
\|f\|_{\hat{X_{\lambda}}^{0,\frac{1}{2},1}}\>
      \|\,|\xi|^{-\frac{1}{2}}g\|_{L^{2}(\SR^{2})}\>
      .
        \label{2.09}
\end{eqnarray}
\end{Lemma}
{\bf Proof.} Combining the identity
\begin{eqnarray*}
&&\tau_{2}+\frac{\xi_{2}^{3}}{4}+\frac{4}{\lambda ^{4}\xi_{2}}
-\left(\tau+\xi^{3}+\frac{1}{\lambda^{4}\xi}\right)
+\left(\tau_{1}+\xi_{1}^{3}+\frac{1}{\lambda^{4}\xi_{1}}\right)\\
&&=\frac{3}{4}\xi_{2}(2\xi-\xi_{2})^{2}
\left[1-\frac{4}{3\lambda^{4}\xi\xi_{1}\xi_{2}^{2}}\right]
\end{eqnarray*}
with a proof similar to  Lemma 2.4 of \cite{LHY},
 we obtain the estimates in Lemma 2.2.

This completes the proof of Lemma 2.2.

\begin{Lemma}\label{Lemma2.3}  The space $\hat{X_{\lambda}}$
 has the following properties:\\
(i) For every $b>1/2$, there exists $C>0$ such that
\begin{eqnarray}
\|f\|_{\hat{X_{\lambda}}}&\leq&
C\|f\|_{\hat{X_{\lambda}}^{-\frac{3}{4},\>b}}.\label{2.010}
\end{eqnarray}
(ii) For $p=2,\frac{4}{3}$, there exists $C>0$ such that
\begin{eqnarray}
&&\|\langle\xi \rangle^{-\frac{3}{4}} f\|_{L_{\xi}^{2}L_{\tau}^{p}}\leq
C\|f\|_{\hat{X_{\lambda}}^{-\frac{3}{4},\>\frac{1}{2}-\epsilon}},\label{2.011}\\
&&\|\langle\xi \rangle^{-\frac{3}{4}} f\|_{\hat{X_{\lambda}}^{-\frac{3}{4},
\>-\frac{1}{2}+\epsilon}}\leq
C\|\langle\xi \rangle^{-\frac{3}{4}}f\|_{L_{\xi}^{2}L_{\tau}^{4}},\label{2.012}\\
&&\|\langle\xi \rangle^{-\frac{3}{4}} f\|_{L_{\xi}^{2}L_{\tau}^{1}}\leq
C\|f\|_{\hat{X_{\lambda}}^{-\frac{3}{4},\frac{1}{2},1}}\label{2.013}.
\end{eqnarray}
\end{Lemma}
\noindent{\bf Proof.} We first prove (i). By using  the Cauchy-Schwartz inequality, since  $b>\frac12,$ we have
\begin{eqnarray}
&&\|u\|_{\hat{X_{\lambda}}^{-\frac{3}{4},\frac{1}{2},1}}^{2}=
\sum_{j}2^{-\frac{3}{2}j}\left(\sum_{k}2^{\frac{k}{2}}
\|u\|_{L_{\xi\tau}^{2}(A_{j}\cap B_{k})}\right)^{2}\nonumber\\
&&\leq \sum_{j}2^{-\frac{3}{2}j}\left(\sum_{k}2^{bk}2^{(\frac{1}{2}-b)k}
\|u\|_{L_{\xi\tau}^{2}(A_{j}\cap B_{k})}\right)^{2}\nonumber\\
&&\leq  \sum_{j}2^{-\frac{3}{2}j}\left(\sum_{k}2^{2bk}
\|u\|_{L_{\xi\tau}^{2}(A_{j}\cap B_{k})}^{2}\right)\left(\sum_{k}2^{(1-2b)k}\right)\nonumber\\
&&\leq \sum_{j}2^{-\frac{3}{2}j}\left(\sum_{k}2^{2bk}\|u\|_{L_{\xi\tau}^{2}(A_{j}\cap B_{k})}^{2}\right)=\|u\|_{\hat{X_{\lambda}}^{-\frac{3}{4},b}}^{2}.\label{2.014}
\end{eqnarray}
Note that
\begin{eqnarray}
\|u\|_{\hat{X_{\lambda}}^{-\frac{3}{4},\frac{1}{2}}}\leq
 C\|u\|_{\hat{X_{\lambda}}^{-\frac{3}{4},b}}.\label{2.015}
\end{eqnarray}
Combining (\ref{2.014}) with (\ref{2.015}),
we have
\begin{eqnarray*}
\|u\|_{\hat{X_{\lambda}}}=\|\mathscr{F}^{-1}u\|_{X_{\lambda}}
=\left\|\mathscr{F}^{-1}[\chi_{D^{c}}
\mathscr{F}u]\right\|_{X_{\lambda}^{-\frac{3}{4}, \frac{1}{2}, 1}}
    +\|\mathscr{F}^{-1}[\chi_{ D}\mathscr{F}u]
    \|_{X_{\lambda}^{-\frac{3}{4}, \frac{1}{2}}}\leq C
\|u\|_{\hat{X_{\lambda}}^{-\frac{3}{4},b}}.
\end{eqnarray*}
The completes the  proof of (i).

Now we prove (ii). The inequality  (\ref{2.013})  is in \cite{BT}.
We only  prove the case $p=\frac{4}{3}$ of (\ref{2.011}),  as the  case $p=2$ of (\ref{2.011})
 is easily checked.
By    the H\"older inequality,  we have
\begin{eqnarray*}
&&\|f\|_{L_{\xi}^{2}L_{\tau}^{\frac{4}{3}}}^{2}=\int_{\SR}
\left(\int_{\SR}|f|^{\frac{4}{3}}d\tau\right)^{\frac{3}{2}}d\xi\nonumber\\&&
=\int_{\SR}\left(\int_{\SR}\langle\sigma^{\lambda}\rangle^{-\frac{2}{3}(1-2\epsilon)}
\langle\sigma^{\lambda}\rangle^{\frac{2}{3}(1-2\epsilon)}
|f|^{\frac{4}{3}}d\tau\right)^{\frac{3}{2}}d\xi
\nonumber\\&&\leq \left[\sup\limits_{\xi \in \SR}
\left(\int_{\SR}\langle\sigma^{\lambda}\rangle
^{-2(1-2\epsilon)}d\tau\right)^{\frac{1}{2}}\right]
\left(\int_{\SR^{2}}\langle\sigma^{\lambda}\rangle^{(1-2\epsilon)}|f|^{2}
d\tau d\xi\right)\nonumber\\
&&\leq C\int_{\SR^{2}}\langle\sigma^{\lambda}\rangle^{(1-2\epsilon)}|f|^{2}d\tau d\xi.
 \end{eqnarray*}
Therefore, from the above  inequality, we have
\begin{eqnarray}
&&\hspace{-1cm}\|\langle\xi \rangle^{-\frac{3}{4}}f\|_{L_{\xi}^{2}L_{\tau}^{\frac{4}{3}}}
=\left[\int_{\SR}\langle\xi \rangle^{-\frac{3}{2}}\left(\int_{\SR}|f|^{\frac{4}{3}}d\tau\right)^{\frac{3}{2}}
d\xi\right]^{\frac12}\nonumber\\&&
\leq C\left(\int_{\SR^{2}}\langle\xi \rangle^{-\frac{3}{2}}
\langle\sigma^{\lambda}\rangle^{(1-2\epsilon)}|f|^{2}d\tau d\xi\right)^{\frac{1}{2}}
=C\left\|\langle\xi\rangle^{-\frac34}f\right\|_{\hat{X_{\lambda}}^{-\frac{3}{4},
\>\frac{1}{2}-\epsilon}}.\label{2.016}
 \end{eqnarray}
By duality, from (\ref{2.016}) we derive that (\ref{2.012}) is valid.
This completes the  proof of (ii).

This completes the proof of Lemma 2.3.

\begin{Lemma}\label{Lemma2.4}  Let $e^{-t(-\partial_{x}^{3}-\lambda^{-4}\partial_{x}^{-1})}v(x,0)$
be the solution to the linear version (discard the nonlinear part) of  equation (\ref{1.03}).
 Then  we have the following estimate
\begin{eqnarray*}
&&\hspace{-1cm}\left\|\psi\left(t\right)
e^{-t(-\partial_{x}^{3}-\lambda^{-4}\partial_{x}^{-1})}v(x,0)\right\|_{X_{\lambda},T}+
\|\psi\left(t\right)e^{-t(-\partial_{x}^{3}-\lambda ^{-4}\partial_{x}^{-1})}
v(x,0)\|_{H_{x}^{-\frac{3}{4}}(\SR)}\nonumber\\&&\leq
C\|v(x,0)\|_{H_{x}^{-\frac{3}{4}}(\SR)}.
\end{eqnarray*}
\end{Lemma}
\noindent{\bf Proof.}  Without loss of generality,  we assume that $T=1$. By   the embedding $X_{\lambda}^{-\frac{3}{4},\frac{1}{2}
+\epsilon}\hookrightarrow X_{\lambda}$ and Lemma 3.1 of \cite{KPV1993}, we have
\begin{eqnarray*}
&&\left\|\psi\left(t\right)e^{-t(-\partial_{x}^{3}-\lambda^{-4}\partial_{x}^{-1})}v(x,0)
\right\|_{X_{\lambda},1}\leq
C\left\|\psi\left(t\right)e^{-t(-\partial_{x}^{3}-\lambda^{-4}\partial_{x}^{-1})}v(x,0)
\right\|_{X_{\lambda}^{-\frac{3}{4},\frac{1}{2}+\epsilon}}\nonumber\\&&\leq C
\|v(x,0)\|_{H_{x}^{-\frac{3}{4}}(\SR)}.
\end{eqnarray*}
Thus,
$
\left\|\psi\left(t\right)e^{-t(-\partial_{x}^{3}-\lambda ^{-4}\partial_{x}^{-1})}v(x,0)
\right\|_{H_{x}^{-\frac{3}{4}}(\SR)}\leq
C\|v(x,0)\|_{H_{x}^{-\frac{3}{4}}(\SR)}.
$

This completes the proof of Lemma 2.4.

\begin{Lemma}\label{Lemma2.5}
The following estimate holds:
\begin{eqnarray*}
&&\left\|\psi\left(t\right)\int_{0}^{t}e^{-(t-s)(-\partial_{x}^{3}-\lambda^{-4}\partial_{x}^{-1})}
         F(s)ds\right\|_{X_{\lambda},T}+\left\|\psi\left(t\right)\int_{0}^{t}
         e^{-(t-s)(-\partial_{x}^{3}-\lambda^{-4}\partial_{x}^{-1})}
           F(s)ds\right\|_{H_{x}^{-\frac{3}{4}}(\SR)}
            \qquad
             \nonumber\\
&&\leq
 C\left\|\mathscr{F}^{-1}
   \Bigl(\left\langle\sigma^{\lambda}\right\rangle^{-1}\mathscr{F}F\Bigr)\right\|_{X_{\lambda}}
   +\left\|\mathscr{F}^{-1}\Bigl(\left\langle\sigma^{\lambda}\right\rangle^{-1}\mathscr{F}F\Bigr)\right\|_{Y}.
\end{eqnarray*}
\end{Lemma}
\noindent{\bf Proof.}  Without loss of generality,  we assume that $T=1$. Note that
\begin{eqnarray*}
\left\|\psi\left(t\right)\int_{0}^{t}
         e^{-(t-s)(-\partial_{x}^{3}-\lambda^{-4}\partial_{x}^{-1})}
           F(s)ds\right\|_{H_{x}^{-\frac34}(\SR)}\leq C\left\|\mathscr{F}^{-1}
           \Bigl(\left\langle\sigma^{\lambda}\right\rangle^{-1}\mathscr{F}F\Bigr)\right\|_{Y}.
\end{eqnarray*}
By  using the Fourier transformation with respect to the space variable, we have that
\begin{eqnarray*}
&&\mathscr{F}_{x}\left[\int_{0}^{t}e^{-(t-s)(-\partial_{x}^{3}-\lambda^{-4}
\partial_{x}^{-1})} F(s)ds\right](\xi)=Ce^{-it(\xi^{3}
    +\frac{1}{\lambda^{4}\xi})}\int_{\SR}\frac{e^{it(\tau+\xi^{3}
    +\frac{1}{\lambda^{4}\xi})}-1}{\tau+\xi^{3}
    +\frac{1}{\lambda^{4}\xi}}\mathscr{F}(\tau,\xi)d\tau\nonumber\\&&=I_{1}+I_{2}+I_{3},
\end{eqnarray*}
where
\begin{eqnarray*}
&&I_{1}=Ce^{-it(\xi^{3}
    +\frac{1}{\lambda^{4}\xi})}\int_{\SR}\frac{e^{it(\tau+\xi^{3}
    +\frac{1}{\lambda^{4}\xi})}-1}{\tau+\xi^{3}
    +\frac{1}{\lambda^{4}\xi}}\left(I_{B_{0}}\mathscr{F}\right)(\tau,\xi)d\tau,\\
 &&I_{2}=C\int_{\SR}\frac{e^{it\tau}}{\tau+\xi^{3}
    +\frac{1}{\lambda^{4}\xi}}\left(I_{B_{>0}}\mathscr{F}\right)(\tau,\xi)d\tau,\\
 &&I_{3}=C\int_{\SR}\frac{e^{-it(\xi^{3}
    +\frac{1}{\lambda^{4}\xi})}}{\tau+\xi^{3}
    +\frac{1}{\lambda^{4}\xi}}\left(I_{B_{>0}}\mathscr{F}\right)(\tau,\xi)d\tau.
\end{eqnarray*}
By using a proof similar to  \cite{GTV} and Lemma 4.1 of \cite{Kis},
we obtain that
\begin{eqnarray*}
&&\left\|I_{1}+I_{2}+I_{3}\right\|_{X_{\lambda}}\leq C\left\|\mathscr{F}^{-1}
   \Bigl(\left\langle\sigma^{\lambda}\right\rangle^{-1}\mathscr{F}F\Bigr)\right\|_{X_{\lambda}}
   +\left\|\mathscr{F}^{-1}\Bigl(\left\langle\sigma^{\lambda}
   \right\rangle^{-1}\mathscr{F}F\Bigr)\right\|_{Y}.
\end{eqnarray*}

This completes the proof of Lemma 2.5.

\begin{Lemma}\label{Lemma2.6}
The following embeddings are true:
\begin{eqnarray*}
X_{\lambda}^{0,\frac{1}{2},1}\hookrightarrow X_{\lambda}^{0,\frac{1}{2}},
X_{\lambda}^{0,\frac{1}{2},1}\hookrightarrow  C(\R;L^{2}(\R)).
\end{eqnarray*}
\end{Lemma}
\noindent {\bf Proof.} The conclusion of Lemma can be found in \cite{Kis}, but the proof is not given. Now we give the explicit proof.
Note that
\begin{eqnarray*}
&&\|u\|_{X_{\lambda}^{0,\frac{1}{2}}}^{2}=\sum_{j}\sum_{k}2^{k}\|\mathscr{F}u\|_{L_{\xi\tau}^{2}(A_{j}\cap B_{k})}^{2},\|u\|_{X_{\lambda}^{0,\frac{1}{2},1}}^{2}=\sum_{j}\left(\sum_{k}2^{\frac{k}{2}}
\|\mathscr{F}u\|_{L_{\xi\tau}^{2}(A_{j}\cap B_{k})}\right)^{2}.
\end{eqnarray*}
To prove $ X_{\lambda}^{0,\frac{1}{2},1}\hookrightarrow X_{\lambda}^{0,\frac{1}{2}}$,
 it suffices to show that
\begin{eqnarray}
\sum_{k}2^{k}\|\mathscr{F}u\|_{L_{\xi\tau}^{2}(A_{j}\cap B_{k})}^{2}\leq
 \left(\sum_{k}2^{\frac{k}{2}}\|\mathscr{F}u\|_{L_{\xi\tau}^{2}(A_{j}\cap B_{k})}\right)^{2}.\label{2.017}
\end{eqnarray}
Let $2^{\frac{k}{2}}\|\mathscr{F}u\|_{L_{\xi\tau}^{2}(A_{j}\cap B_{k})}=a_{k}$. Then, we have
\begin{eqnarray}
\left[\sum_{k}a_{k}^{2}\right]^{\frac{1}{2}}\leq \sum_{k}a_{k},\label{2.018}
\end{eqnarray}
which yields that (\ref{2.017})  is valid. Thus, we have that
$X_{\lambda}^{0,\frac{1}{2},1}\hookrightarrow X_{\lambda}^{0,\frac{1}{2}}.$

\noindent Now we prove that $X_{\lambda}^{0,\frac{1}{2},1}\hookrightarrow  C(\R;L^{2}(\R)).$
\noindent By using (\ref{2.013}), we have that
\begin{eqnarray}
\|u(t)\|_{L^{2}}\leq C\left\|\int_{\SR}e^{it \tau}
\mathscr{F}u(\xi,\tau)d\tau\right\|_{L_{\xi}^{2}}\leq
C\|\mathscr{F}u(\xi,\tau)\|_{L_{\xi}^{2}L_{\tau}^{1}}
\leq C\|u\|_{X_{\lambda}^{0,\frac{1}{2},1}}.\label{2.019}
\end{eqnarray}
 For $\forall \epsilon>0,$  from (\ref{2.013}), we know that
  there exists  sufficiently large $A>0$ such that
 \begin{eqnarray}
 \left\|\int_{|\tau|\geq A}|\mathscr{F}u(\xi,\tau)|d\tau\right\|_{L_{\xi}^{2}} <\frac{\epsilon}{4}.\label{2.020}
 \end{eqnarray}
For this  $\epsilon>0,$ $\exists \delta$, which satisfies
 $0<\delta<\frac{\epsilon}{2A\|u\|_{X_{\lambda}^{0,\frac{1}{2},1}}},$  and when
$|t_{1}-t_{2}|<\delta,$ we infer that
\begin{eqnarray}
&&\|u(t_{1})-u(t_{2})\|_{L^{2}}\leq\left\|\int_{|\tau|\geq A}
(e^{it_{1} \tau}-e^{it_{2} \tau})
\mathscr{F}u(\xi,\tau)d\tau\right\|_{L_{\xi}^{2}}\nonumber\\&&\qquad
+\left\|\int_{|\tau|< A}(e^{it_{1} \tau}-e^{it_{2} \tau})
\mathscr{F}u(\xi,\tau)d\tau\right\|_{L_{\xi}^{2}}\leq
\left\|\int_{|\tau|\geq A}|e^{it_{1} \tau}-e^{it_{2} \tau}|
\mathscr{F}u(\xi,\tau)d\tau\right\|_{L_{\xi}^{2}}\nonumber\\&&\qquad
+\left\|\int_{|\tau|< A}|e^{it_{1} \tau}-e^{it_{2} \tau}|
\mathscr{F}u(\xi,\tau)d\tau\right\|_{L_{\xi}^{2}}\nonumber\\&&\leq
2\left\|\int_{|\tau|\geq A}|\mathscr{F}u(\xi,\tau)|d\tau\right\|_{L_{\xi}^{2}}+
\left\|\int_{|\tau|< A}|t_{1}-t_{2}||\tau|\mathscr{F}u(\xi,\tau)d\tau\right\|_{L_{\xi}^{2}}\nonumber\\
&&\leq
\frac{\epsilon}{2}+|t_{1}-t_{2}|A\left\|\int_{|\tau|< A}\mathscr{F}u(\xi,\tau)d\tau\right\|_{L_{\xi}^{2}}\nonumber\\
&&\leq \frac{\epsilon}{2}+A\delta\left\|\int_{|\tau|< A}|\mathscr{F}u(\xi,\tau)|d\tau\right\|_{L_{\xi}^{2}}\nonumber\\
&&\leq \frac{\epsilon}{2}+A\delta\|u\|_{L_{\xi}^{2}L_{\tau}^{1}}\nonumber\\
&&\leq \frac{\epsilon}{2}+A\delta\|u\|_{X_{\lambda}^{0,\frac{1}{2},1}}\leq \frac{\epsilon}{2}+\frac{\epsilon}{2}=\epsilon.\label{2.021}
\end{eqnarray}
Thus, putting  (\ref{2.019})-(\ref{2.021})  together,
we conclude  that $X_{\lambda}^{0,\frac{1}{2},1}\hookrightarrow  C(\R;L^{2}(\R)).$

This ends the proof of   Lemma 2.6.

\begin{Lemma}\label{Lemma2.7}
Take  $\epsilon$ with  $0<\epsilon<\frac{1}{10^{8}}$ and let $\mathscr{F}(P^{a}f)(\xi)=\chi_{\{|\xi|\geq a\}}(\xi)
\mathscr{F}f(\xi)$ with $a\geq4$. Then, we have
\begin{eqnarray}
&&\|u\|_{L_{xt}^{6}}\leq C\|u\|_{X_{\lambda}^{0,\frac{1}{2}+\epsilon}},\label{2.022}\\
&&\left\|D_{x}^{\frac{1}{6}}P^{a}u\right\|_{L_{xt}^{6}}\leq C
\|u\|_{X_{\lambda}^{0,\frac{1}{2}+\epsilon}},\label{2.023}\\
&&\|u\|_{L_{xt}^{4}}\leq C\|u\|_{X_{\lambda}^{0,\frac{3}{4}
\left(\frac{1}{2}+\epsilon\right)}},\label{2.024}\\
&&\left\|D_{x}^{\frac{1}{8}}P^{a}u\right\|_{L_{xt}^{4}}\leq
 C\|u\|_{X_{\lambda}^{0,\frac{3}{4}(\frac{1}{2}+\epsilon)}}\label{2.025}.
\end{eqnarray}
\end{Lemma}
\noindent{\bf Proof.}
For the proof of (\ref{2.022}) and (\ref{2.024}), we refer the readers
 to (2.27) and (2.21) of \cite{GL} and (2.3) of Lemma 2.1 in \cite{GH}.
Interpolating (\ref{2.022}), (\ref{2.023}) with $\|u\|_{L_{xt}^{2}}
=\|u\|_{X_{0,0}}$ yields that (\ref{2.024})-(\ref{2.025}) are valid.

This completes the proof of Lemma 2.7.

\begin{Lemma}\label{Lemma2.8}
Let $T\in (0,1)$ and $-\frac{1}{2}<b^{\prime}\leq0\leq b\leq b^{\prime}+1$. Then the following estimates hold for  $s\in \R$
\begin{eqnarray}
&&\|\eta_{T}(t)S^{\lambda}(t)\phi\|_{X_{\lambda}^{s,b}(\SR^{2})}\leq
 CT^{\frac{1}{2}-b}\|\phi\|_{H^{s}(\SR)},\nonumber\\
&&\left\|\eta_{T}(t)\int_{0}^{t}S^{\lambda}(t-\tau)F(\tau)d\tau\right\|_{X_{\lambda}^{s,b}(\SR^{2})}
\leq CT^{1+b^{\prime}-b}\|F\|_{X_{\lambda}^{s,b^{\prime}}(\SR^{2})}.\nonumber
\end{eqnarray}
\end{Lemma}

For the proof of Lemma 2.8, we refer the readers to
 \cite{Bourgain-GAFA93,KPV1993,Grunrock}.

\begin{Lemma}\label{Lemma2.9}
We define $I^{s}$  as follows
\begin{eqnarray*}
\mathscr{F}I^{s}_{-}\left(u_{1},u_{2}\right)=C\int_{\tiny\begin{array}{cc}\xi
=\sum\limits_{j=1}^{2}\xi_j\\ \tau =\sum\limits_{j=1}^{2}\tau_j\end{array}}
\left|\left[\phi^{\lambda}(\xi_{1})\right]^{\prime}-\left[\phi^{\lambda}
(\xi_{2})\right]^{\prime}\right|^{s}\mathscr{F}u_{1}(\xi_{1},\tau_{1})
\mathscr{F}u_{2}(\xi_{2},\tau_{2}) d\xi_{1}d\tau_{1}.
\end{eqnarray*}
Let $0<s<\frac{1}{2}$ and $b=\frac{1}{2}+\epsilon$.  Then,
\begin{eqnarray*}
\left\|I^{s}_{-}\left(u_{1},u_{2}\right)\right\|_{L_{xt}^{2}}\leq C\prod_{j=1}^{2}
\|u_{j}\|_{X_{\lambda}^{0,\frac{3+2s}{4}b}}.
\end{eqnarray*}
\end{Lemma}
\noindent{\bf Proof.} From Lemma 2.5 of \cite{LW},  we  conclude that
\begin{eqnarray}
\left\|I^{\frac{1}{2}}_{-}\left(u_{1},u_{2}\right)\right\|_{L_{xt}^{2}}\leq C\prod_{j=1}^{2}
\|u_{j}\|_{X_{\lambda}^{0,\frac{1}{2}+\epsilon}}.\label{2.026}
\end{eqnarray}
Let $F_{j}(\xi_{j},\tau_{j})=\langle \sigma_{j}^{\lambda}\rangle^{\frac{3+2s}{4}b}
\mathscr{F}u_{j}(\xi_{j},\tau_{j})(j=1,2).$
By   the Plancherel identity,  we know that in order to obtain (\ref{2.026}), it suffices
to prove that
\begin{eqnarray}
\left\|\int_{\tiny\begin{array}{cc}\xi=\sum\limits_{j=1}^{2}\xi_j\\ \tau =\sum\limits_{j=1}^{2}\tau_j\end{array}}\left|\left[\phi^{\lambda}(\xi_{1})\right]^{\prime}-\left[\phi^{\lambda}
(\xi_{2})\right]^{\prime}\right|^{s}\prod_{j=1}^{2}\left(\frac{F_{j}}{\langle \sigma_{j}^{\lambda}\rangle^{\frac{3+2s}{4}b}}\right)d\xi_{1}d\tau_{1}\right\|_{L_{\xi\tau}^{2}}
\leq C\prod_{j=1}^{2}\|F_{j}\|_{L_{\xi\tau}^{2}}.\label{2.027}
\end{eqnarray}
Define $b_{1}=\frac{3+2s}{4}b.$ By using the Young inequality and the fact that $0<s<\frac{1}{2}$,
we derive that
\begin{eqnarray}
&&\left|\left[\phi^{\lambda}(\xi_{1})\right]^{\prime}-\left[\phi^{\lambda}
(\xi_{2})\right]^{\prime}\right|^{s}\prod_{j=1}^{2}\langle\sigma_{j}^{\lambda}\rangle^{-b_{1}}\nonumber\\
&&\leq \left|\left[\phi^{\lambda}(\xi_{1})\right]^{\prime}-\left[\phi^{\lambda}
(\xi_{2})\right]^{\prime}\right|^{s}\left(\prod_{j=1}^{2}\langle\sigma_{j}^{\lambda}\rangle^{-2bs}\right)
\left(\prod_{j=1}^{2}\langle\sigma_{j}^{\lambda}\rangle^{-(b_{1}-2bs)}\right)\nonumber\\
&&\leq 2s\left|\left[\phi^{\lambda}(\xi_{1})\right]^{\prime}-\left[\phi^{\lambda}
(\xi_{2})\right]^{\prime}\right|^{\frac{1}{2}}\left(\prod_{j=1}^{2}\langle\sigma_{j}^{\lambda}\rangle^{-b}\right)+(1-2s)
\left(\prod_{j=1}^{2}\langle\sigma_{j}^{\lambda}\rangle^{-\frac{3}{4}b}\right).\label{2.028}
\end{eqnarray}
Combining (\ref{2.024}), (\ref{2.026}) with (\ref{2.028}),  we have that
\begin{eqnarray}
&&\left\|\int_{\tiny\begin{array}{cc}\xi=\sum\limits_{j=1}^{2}\xi_j\\ \tau =\sum\limits_{j=1}^{2}\tau_j\end{array}}\left|\left[\phi^{\lambda}(\xi_{1})\right]^{\prime}-\left[\phi^{\lambda}
(\xi_{2})\right]^{\prime}\right|^{s}\prod_{j=1}^{2}\left(\frac{F_{j}}{\langle \sigma_{j}^{\lambda}\rangle^{\frac{3+2s}{4}b}}\right)d\xi_{1}d\tau_{1}\right\|_{L_{\xi\tau}^{2}}\nonumber\\&&\leq
\left\|\int_{\tiny\begin{array}{cc}\xi=\sum\limits_{j=1}^{2}\xi_j\\ \tau =\sum\limits_{j=1}^{2}\tau_j\end{array}}\left|\left[\phi^{\lambda}(\xi_{1})\right]^{\prime}-\left[\phi^{\lambda}
(\xi_{2})\right]^{\prime}\right|^{\frac{1}{2}}\prod_{j=1}^{2}\left(\frac{F_{j}}{\langle \sigma_{j}^{\lambda}\rangle^{b}}\right)d\xi_{1}d\tau_{1}\right\|_{L_{\xi\tau}^{2}}\nonumber\\&&\qquad\qquad
+\left\|\int_{\tiny\begin{array}{cc}\xi=\sum\limits_{j=1}^{2}\xi_j\\ \tau
 =\sum\limits_{j=1}^{2}\tau_j\end{array}}\prod_{j=1}^{2}\left(\frac{F_{j}}{\langle \sigma_{j}^{\lambda}\rangle^{\frac{3}{4}b}}\right)d\xi_{1}d\tau_{1}\right\|_{L_{\xi\tau}^{2}}\nonumber\\
 &&\leq C\prod_{j=1}^{2}\|F_{j}\|_{L_{\xi\tau}^{2}}
.\label{2.029}
\end{eqnarray}

This completes the proof of  Lemma 2.9.

\bigskip
\bigskip

\noindent{\large\bf 3. Bilinear estimates for $s=-\frac34$}

\setcounter{equation}{0}

 \setcounter{Theorem}{0}

\setcounter{Lemma}{0}

 \setcounter{section}{3}

In this section,  we establish the  bilinear estimates  in the case of  $s=-\frac34$.
These will be needed in proving Theorem 1.1 (well-posedness).
More precisely, by using Lemmas 2.1-2.3, 2.6,  we first establish Lemma 3.1 and then apply Lemma 3.1 to
establish Lemma 3.2.

\begin{Lemma}\label{Lemma3.1}  Assume that  $f, g \in \mathscr{S}'(\R^{2})$, with
$\supp f\subset A_{j_{1}}$ and $\supp g\subset A_{j_{2}}$. Then
\begin{eqnarray}
     \left\|\chi_{A_{j}}\left\langle\sigma^{\lambda}\right\rangle^{-1}
      \xi f*g\right\|_{\hat{X_{\lambda}}}
&\leq&  C(j,j_{1},j_{2})\|f\|_{\hat{X_{\lambda}}}\|g\|_{\hat{X_{\lambda}}},
        \label{3.01}\\
     \left\|\chi_{A_{j}}\langle\xi\rangle^{-3/4}
      \left\langle\sigma^{\lambda}\right\rangle^{-1}
       \xi f*g\right\|_{L_{\xi}^{2}L_{\tau}^{1}}
&\leq& C(j,j_{1},j_{2})\|f\|_{\hat{X_{\lambda}}}\|g\|_{\hat{X_{\lambda}}}
         \label{3.02}
\end{eqnarray}
for $j\geq0$ in the following cases.

(i)  At least two of $j,j_{1},j_{2}$ are less than 40 and $C(j,j_{1},j_{2})\sim 1$.

(ii) $j_{1},j_{2}\geq40$,  $|j_{1}-j_{2}|\leq 10$,   $0<j<j_{1}-9$
 and $C(j,j_{1},j_{2})\sim 2^{-\frac{3}{8} j}$.

(iii) $j,j_{1}\geq 40,$ $|j-j_{1}|\leq 10$, $0<j_{2}<j-10$  and
  $C(j,j_{1},j_{2})\sim2^{-\frac{1}{4}(j-j_{2})}.$

(iv) $j,j_{2}\geq 40,$ $|j-j_{2}|\leq 10$, $0<j_{1}<j-10$ and
 $C(j,j_{1},j_{2})\sim2^{-\frac{1}{4}(j-j_{1})}.$

(v) $j,j_{1},j_{2}\geq 40$, $|j-j_{1}|\leq 10$, $|j-j_{2}|\leq 10$
 and $C(j,j_{1},j_{2})\sim 1.$

(vi) $j_{1},j_{2}\geq40$,   $j=0$ and $C(j,j_{1},j_{2})\sim 1.$

(vii)  $j,j_{1}\geq 40$,   $j_{2}=0$ and $C(j,j_{1},j_{2})\sim 1.$

(viii)  $j,j_{2}\geq 40$,  $j_{1}=0$ and $C(j,j_{1},j_{2})\sim 1.$

\end{Lemma}

{\bf Proof.} (i)  This case can be proved similarly to (i) of \cite{LHY}
 with the aid of Lemma 2.3.

\noindent(ii) In this case,
 we claim that
\begin{eqnarray}
\left|\sigma^{\lambda}-\sigma_{1}^{\lambda}
-\sigma_{2}^{\lambda}\right|
=\left|3\xi\xi_{1}\xi_{2}-\frac{\xi_{1}^{2}+\xi_{1}\xi_{2}
+\xi_{2}^{2}}{\lambda ^{4}\xi\xi_{1}\xi_{2}}\right|\geq \frac{|\xi\xi_{1}\xi_{2}|}{4}.\label{3.03}
\end{eqnarray}
Indeed, if otherwise
\begin{eqnarray}
\left|\sigma^{\lambda}-\sigma_{1}^{\lambda}
-\sigma_{2}^{\lambda}\right|
=\left|3\xi\xi_{1}\xi_{2}-\frac{\xi_{1}^{2}+\xi_{1}\xi_{2}
+\xi_{2}^{2}}{\lambda ^{4}\xi\xi_{1}\xi_{2}}\right|< \frac{|\xi\xi_{1}\xi_{2}|}{4}\label{3.04}
\end{eqnarray}
is valid, then
 we obtain that
\begin{eqnarray}
3|\xi\xi_{1}\xi_{2}|-\frac{\xi_{1}^{2}+\xi_{1}\xi_{2}+\xi_{2}^{2}}
{|\lambda ^{4}\xi\xi_{1}\xi_{2}|}< \frac{|\xi\xi_{1}\xi_{2}|}{4}.   \label{3.05}
\end{eqnarray}
But we now verify that this will lead to a contradiction.
By  (\ref{3.05}) and the  fact  $\lambda \geq 1$, we have
\begin{eqnarray}
\frac{11}{4}\xi^{2}\xi_{1}^{2}\xi_{2}^{2}\leq \frac{11}{4}\lambda ^{4}\xi^{2}
\xi_{1}^{2}\xi_{2}^{2}\leq \xi_{1}^{2}+\xi_{1}\xi_{2}+\xi_{2}^{2}.  \label{3.06}
\end{eqnarray}
From (\ref{3.06}), we have
\begin{eqnarray}
\frac{11}{4}\xi^{2}\xi_{2}^{2}\leq 1+\frac{\xi_{1}}{\xi_{2}}
+\left|\frac{\xi_{1}}{\xi_{2}}\right|^{2}\leq 7,\label{3.07}
\end{eqnarray}
which contradicts with $|\xi|\geq 1$ and $|\xi_{2}|\geq 100.$
Thus,  our claim (\ref{3.03}) is true. It  is easily checked that
 Lemmas 2.1-2.2 are true. We restrict  $f$ to $B_{k_{1}}$ and $g$ to $B_{k_{2}}$.
From  (\ref{3.03}), we have
\begin{eqnarray}
2^{k_{max}}:=2^{{\rm max}\{k,\>k_{1},\>k_{2}\}}\geq C2^{j+2j_{1}}\label{3.08}.
\end{eqnarray}
We consider
\begin{eqnarray*}
 2^{k_{1}}=2^{{\rm max}\{k,\>k_{1},\>k_{2}\}},2^{k_{2}}=2^{{\rm max}\{k,\>k_{1},\>k_{2}\}},2^{k}=2^{{\rm max}\{k,\>k_{1},\>k_{2}\}},
\end{eqnarray*}
respectively.

\noindent When $2^{k_{1}}\geq C2^{j+2j_{1}}$ which leads to that $2^{-\frac{3}{8}(j+2j_{1})}2^{\frac{3}{8}k_{1}}
2^{\frac{1}{8}(k_{1}-k)}\geq C.$  By  Lemma 2.2 and a proof similar to
$2^{k_{1}}\geq C2^{j+2j_{1}}$  of (ii) in \cite{LHY},  we have
\begin{eqnarray}
&&\hspace{-0.5cm}\left\|I_{A_{j}}\left\langle\sigma^{\lambda}\right\rangle^{-1}\xi f*g\right\|_{\hat{X_{\lambda}}}
\leq C2^{-\frac{3j}{8}}\|f\|_{\hat{X_{\lambda}}}
\|g\|_{\hat{X_{\lambda}}}\label{3.09}.
\end{eqnarray}
When $2^{k_{2}}\geq C2^{j+2j_{1}}$,  this case can be treated similarly
to the case $2^{k_{1}}\geq C2^{j+2j_{1}}$.

\noindent When $2^{k}=2^{{\rm max}\{k,\>k_{1},\>k_{2}\}}$, we consider $2^{k}\gg 2^{{\rm max}\{k_{1},\>k_{2}\}}$ and
$2^{k}\sim 2^{{\rm max}\{k_{1},\>k_{2}\}}$, respectively.

\noindent When $2^{k}\gg 2^{{\rm max}\{k_{1},\>k_{2}\}}$,  we have that $2^{k}\sim 2^{j+2j_{1}}$ which leads to   $2^{\frac{j}{4}}2^{-\frac{k}{2}}
\leq C2^{-\frac{3j}{4}}2^{-j_{1}}2^{\frac{j}{2}}.$
By   Lemma 2.1 and a proof similar to  (ii) of \cite{LHY},  we have
\begin{eqnarray}
&&\left\|I_{A_{j}}\left\langle\sigma^{\lambda}\right\rangle^{-1}\xi
f*g\right\|_{\hat{X_{\lambda}}}\leq C2^{-\frac{3}{8}j}\|f\|_{\hat{X_{\lambda}}}
\|g\|_{\hat{X_{\lambda}}}
\label{3.010}.
\end{eqnarray}
When $2^{k}\sim 2^{{\rm max}\{k_{1},\>k_{2}\}}$,  we have that $2^{k}\sim 2^{k_{1}}$ or $2^{k}\sim 2^{k_{2}}$, respectively.

\noindent When $2^{k}\sim 2^{k_{1}}$, this case can be proved similarly to case $2^{k_{1}}=2^{{\rm max}\{k,\>k_{1},\>k_{2}\}}$.

\noindent When $2^{k}\sim 2^{k_{2}}$, this case can be proved similarly to case $2^{k_{2}}=2^{{\rm max}\{k,\>k_{1},\>k_{2}\}}$.

With the aid of the Cauchy-Schwartz inequality with respect to $\tau$, and
combining the
case $2^{k}\geq C2^{j+2j_{1}}$ with $2^{k_{1}}\geq C2^{j+2j_{1}}$
and $2^{k_{2}}\geq C2^{j+2j_{1}}$, we have
\begin{eqnarray}
&&\left\|I_{A_{j}}\left\langle\sigma^{\lambda}\right\rangle^{-1}
\xi \langle \xi\rangle ^{-\frac{3}{4}}f*g\right\|_{L_{\xi}^{2}L_{\tau}^{1}}\nonumber\\&&
\leq C2^{\frac{j}{4}}\sum_{k\geq0}2^{-\frac{k}{2}}\left\| f*g\right\|
_{L_{\xi}^{2}L_{\tau}^{2}}\nonumber\\&&
\leq C\left[2^{-\frac{3}{4}j}2^{-\frac{3}{2}j_{1}}
+2^{-\frac{1}{8}j}2^{-\frac{7}{4}j_{1}}\right]\|f\|_{\hat{X_{\lambda}}^{0,\frac{1}{2},1}}
\|g\|_{\hat{X_{\lambda}}^{0,\frac{1}{2},1}}\nonumber\\&&
\leq C2^{-\frac{3}{8}j}
\|f\|_{\hat{X_{\lambda}}^{-\frac{3}{4},\frac{1}{2},1}}\|g\|_{\hat{X_{\lambda}}
^{-\frac{3}{4},\frac{1}{2},1}}\nonumber\\&&\leq C2^{-\frac{3}{8}j}\|f\|_{\hat{X_{\lambda}}}
\|g\|_{\hat{X_{\lambda}}}\label{3.011}.
\end{eqnarray}
(iii)
By a proof similar to (\ref{3.03}),
 we have  $2^{k_{\rm max}}\geq C|\xi
\xi_{1}\xi_{2}|\geq C2^{2j+j_{2}}$.
By using Lemmas 2.1-2.2 and a proof similar to (iii) of \cite{LHY}, we have
\begin{eqnarray}
     \left\|\chi_{A_{j}}\left\langle\sigma^{\lambda}\right\rangle^{-1}
      \xi f*g\right\|_{\hat{X_{\lambda}}}
&\leq&  C2^{-\frac{1}{4}(j-j_{2})}\|f\|_{\hat{X_{\lambda}}}\|g\|_{\hat{X_{\lambda}}}\label{3.012}.
\end{eqnarray}
(iv) This case can be proved similarly to the case (iii).

\noindent(v) By using a proof similar to (\ref{3.03}),
we have that $2^{k_{\rm max}}\geq C|\xi
\xi_{1}\xi_{2}|\geq C2^{3j}\sim 2^{3j_{1}}\sim 2^{3j_{2}}.$
The left hand side of (\ref{3.01})-(\ref{3.02}) can be controlled by
\begin{eqnarray}
C2^{\frac{j}{4}}\sum_{k\geq0}2^{-\frac{k}{2}}\left\| f*g\right\|_{L_{\xi}^{2}L_{\tau}^{2}}\label{3.013}.
\end{eqnarray}
When $2^{k}=2^{k_{\rm max}}\geq C2^{3j}$,  we use   (\ref{2.03})  to estimate (\ref{3.013})
 by
\begin{eqnarray}
C2^{-\frac{3}{2}j}\|f\|_{\hat{X_{\lambda}}^{0,\frac{1}{2},1}}\|g\|_{\hat{X_{\lambda}}^{0,\frac{1}{2},1}}
\leq C\|f\|_{\hat{X_{\lambda}}^{-\frac{3}{4},\frac{1}{2},1}}
\|g\|_{\hat{X_{\lambda}}^{-\frac{3}{4},\frac{1}{2},1}}\leq C\|f\|_{\hat{X_{\lambda}}}
\|g\|_{\hat{X_{\lambda}}}\label{3.014}.
\end{eqnarray}
When $2^{k_{1}}= 2^{k_{\rm max}}\geq C2^{3j}$, by using (\ref{2.08}) and the fact
that $2^{\frac{k_{1}}{2}}2^{-\frac{3j}{2}}\geq C$ as well as $X_{\lambda}^{0,\frac{1}{2},1}\hookrightarrow X_{\lambda}^{0,\frac{1}{2}}$
of Lemma 2.6, (\ref{3.013})  can be controlled by
\begin{eqnarray}
&&C2^{-\frac{5j}{4}}2^{\frac{k_{1}}{2}}\sum_{k\geq0}2^{-\frac{k}{2}}\left\| f*g\right\|_{L_{\xi}^{2}L_{\tau}^{2}}\nonumber\\&&
\leq C2^{-\frac{3}{2}j}\sum_{k\geq 0}2^{-\frac{k}{4}}\|f\|_{\hat{X_{\lambda}}^{0,\frac{1}{2}}}
\|g\|_{\hat{X_{\lambda}}^{0,\frac{1}{2},1}}\nonumber\\&&
\leq C\|f\|_{\hat{X_{\lambda}}^{-\frac{3}{4},\frac{1}{2},1}}\|g\|_{\hat{X_{\lambda}}^{-\frac{3}{4},\frac{1}{2},1}}\nonumber\\&&
\leq C\|f\|_{\hat{X_{\lambda}}}
\|g\|_{\hat{X_{\lambda}}}\label{3.015}.
\end{eqnarray}
When $2^{k_{2}}= 2^{k_{\rm max}}\geq C2^{3j}$, this case can be proved similarly to
case $2^{k_{1}}\sim 2^{k_{\rm max}}\geq C2^{3j}$.

\noindent (vi)
We consider
\begin{eqnarray*}
\left|3\xi\xi_{1}\xi_{2}-\frac{\xi_{1}^{2}+\xi_{1}\xi_{2}+\xi_{2}^{2}}{\lambda^{4}\xi\xi_{1}\xi_{2}}
\right|<\frac{|\xi\xi_{1}\xi_{2}|}{4}, \;\;\;    \left|3\xi\xi_{1}\xi_{2}
-\frac{\xi_{1}^{2}+\xi_{1}\xi_{2}+\xi_{2}^{2}}{\lambda^{4}\xi\xi_{1}\xi_{2}}
\right|\geq\frac{|\xi\xi_{1}\xi_{2}|}{4}.
\end{eqnarray*}
(vi-1) When $\left|3\xi\xi_{1}\xi_{2}-\frac{\xi_{1}^{2}+\xi_{1}\xi_{2}+\xi_{2}^{2}}{\lambda^{4}\xi\xi_{1}\xi_{2}}
\right|<\frac{|\xi\xi_{1}\xi_{2}|}{4}$,
 we have
$|\xi|\sim \lambda^{-2}|\xi_{1}|^{-1}\sim\lambda ^{-2}2^{-j_{1}}.$
By using the H\"older inequality  with respect to $\xi$, the Young inequality with respect
 to $\xi,\tau$ and (\ref{2.013}),  and noting that
$\lambda \geq 1,$ we have
\begin{eqnarray}
&&\left\|I_{A_{j}}\langle\xi\rangle^{-3/4}\left\langle\sigma^{\lambda}\right
\rangle^{-1}\xi f*g\right\|_{L_{\xi}^{2}L_{\tau}^{1}}\leq C\lambda ^{-3}2^{-\frac{3j_{1}}{2}}
\|f*g\|_{L_{\xi}^{\infty}L_{\tau}^{1}}\nonumber\\&&\leq C2^{-\frac{3j_{1}}{2}}\|f\|_{L_{\xi}^{2}
L_{\tau}^{1}}\|g\|_{L_{\xi}^{2}L_{\tau}^{1}}\leq C
\|f\|_{\hat{X_{\lambda}}^{-\frac{3}{4},\frac{1}{2},1}}\|g\|_{\hat{X_{\lambda}}^{-\frac{3}{4},\frac{1}{2},1}}\nonumber\\&&
\leq C\|f\|_{\hat{X_{\lambda}}}
\|g\|_{\hat{X_{\lambda}}}
.\label{3.016}
\end{eqnarray}
When $(\tau,\xi)\in D^{\prime},$ from the definition of $D^{\prime},$ we see  that $|\tau|\geq |\xi|^{-3}
\geq C\lambda ^{6}2^{3j_{1}}$ which implies that $\langle\sigma^{\lambda}\rangle\geq C\lambda ^{6}2^{3j_{1}}$.
By the H\"older inequality  with respect to  $\xi$, the Young
 inequality with respect to $\xi,\tau$ and (\ref{2.011}), (\ref{2.013}), we have
\begin{eqnarray}
&&\hspace{-1cm}\left\|I_{A_{j}}\langle\xi\rangle^{-3/4}\left\langle\sigma^{\lambda}\right
\rangle^{-\frac{1}{2}}\xi f*g\right\|_{L_{\xi}^{2}L_{\tau}^{2}}\nonumber\\&&
\hspace{-1cm}\leq C\lambda^{-5}2^{-\frac{5j_{1}}{2}}
\|f*g\|_{L_{\xi}^{2}L_{\tau}^{2}}\leq C 2^{-3j_{1}}
\|f*g\|_{L_{\xi}^{\infty}L_{\tau}^{2}}\nonumber\\&&\hspace{-1cm}
\leq C\lambda^{-3}2^{-3j_{1}}
\|f\|_{L_{\xi}^{2}L_{\tau}^{2}}\|g\|_{L_{\xi}^{2}L_{\tau}^{1}}
\leq C2^{-\frac{3j_{1}}{2}}
\|f\|_{\hat{X_{\lambda}}^{-\frac{3}{4},\frac{1}{2},1}}
\|g\|_{\hat{X_{\lambda}}^{-\frac{3}{4},\frac{1}{2},1}}\nonumber\\&&\hspace{-1cm}
\leq C\|f\|_{\hat{X_{\lambda}}}
\|g\|_{\hat{X_{\lambda}}}.\label{3.017}
\end{eqnarray}
When $(\tau,\xi)\in D_{1}\cup D_{2}\cup D_{3}$,
 by using the H\"older inequality  with respect to $\xi$ and the Young
inequality and  (\ref{2.011}), (\ref{2.013}) as  well as  $|\xi|\sim \lambda ^{-2}2^{-j_{1}}$,  we have
\begin{eqnarray}
&&\|\xi f*g\|_{\hat{X_{\lambda}}^{-\frac{3}{4},-\frac{1}{2},1}}\leq
C\sum_{k\geq 0}2^{-\frac{k}{2}}\| \xi f*g\|_{L_{\xi\tau}^{2}}
\nonumber\\&&
\leq C\|\xi\|_{L^{2}}\|  f*g\|_{L_{\xi}^{\infty}L_{\tau}^{2}}
\leq C2^{-\frac{3j_{1}}{2}}\|  f*g\|_{L_{\xi}^{\infty}L_{\tau}^{2}}
\leq C2^{-\frac{3j_{1}}{2}}\|f\|_{L_{\xi}^{2}L_{\tau}^{1}}
\|g\|_{L_{\xi}^{2}L_{\tau}^{2}}\nonumber\\
&&\leq\|f\|_{\hat{X_{\lambda}}^{-\frac{3}{4},\frac{1}{2},1}}
\|g\|_{\hat{X_{\lambda}}^{-\frac{3}{4},\frac{1}{2},1}}\leq C\|f\|_{\hat{X_{\lambda}}}
\|g\|_{\hat{X_{\lambda}}}.\label{3.018}
\end{eqnarray}
\noindent (vi-2) When $\left|3\xi\xi_{1}\xi_{2}-\frac{\xi_{1}^{2}+\xi_{1}\xi_{2}+\xi_{2}^{2}}
{\lambda^{4}\xi\xi_{1}\xi_{2}}\right|\geq\frac{|\xi\xi_{1}\xi_{2}|}{4}$, we have that
$2^{k_{\rm max}}=2^{{\rm max}\left\{k,k_{1},k_{2}\right\}}\geq C2^{2j_{1}}|\xi|.$

\noindent In this case, we consider three subcases
\begin{eqnarray*}
(a):2^{k_{1}}=2^{{\rm max}};   \;\;  (b):2^{k_{2}}=2^{{\rm max}};  \;\;   (c):2^{k}= 2^{{\rm max}},
\end{eqnarray*}
 respectively.

\noindent Now we consider  the subcase $(a):  2^{k_{1}}= 2^{k_{\rm max}}$,  and the subcase
 $(b):   2^{k_{2}}= 2^{k_{\rm max}}$. In these subcases, we consider two situations
 \begin{eqnarray*}
 \left|1-\frac{4}{3\lambda^{4}\xi\xi_{1}\xi_{2}^{2}}\right|>  \frac{1}{2}, \;\; \;
 \left|1-\frac{4}{3\lambda^{4}\xi\xi_{1}\xi_{2}^{2}}\right|\leq  \frac{1}{2}.
 \end{eqnarray*}

\noindent When $\left|1-\frac{4}{3\lambda^{4}\xi\xi_{1}\xi_{2}^{2}}\right|>  \frac{1}{2}$, case
 $(a):2^{k_{1}}= 2^{k_{\rm max}}$  can be proved
similarly to case $2^{k_{1}}= 2^{k_{\rm max}}$ and case $2^{k_{2}}= 2^{k_{\rm max}}$
 of (ii) of \cite{LHY}.

\noindent When $\left|1-\frac{4}{3\lambda^{4}\xi\xi_{1}\xi_{2}^{2}}\right|>  \frac{1}{2}$,
   $(b):2^{k_{2}}= 2^{k_{\rm max}}$  can be proved
similarly to case $2^{k_{2}}= 2^{k_{\rm max}}$ of (ii) of \cite{LHY}.

\noindent When  $\left|1-\frac{4}{3\lambda^{4}\xi\xi_{1}\xi_{2}^{2}}\right|\leq  \frac{1}{2}$,
 we have that
\begin{eqnarray}
\left|1+\frac{4}{3\lambda^{4}\xi\xi_{1}^{2}\xi_{2}}\right|
=\left|1-\frac{4}{3\lambda^{4}|\xi\xi_{1}\xi_{2}^{2}|}\right|\leq  \frac{1}{2}\label{3.019}.
\end{eqnarray}
From (\ref{3.019}), we have that
\begin{eqnarray}
  \frac{1}{2}\leq \frac{4}{3\lambda^{4}|\xi\xi_{1}\xi_{2}^{2}|}  \leq \frac{3}{2}.\label{3.020}
 \end{eqnarray}
  From   (\ref{3.020}), we have that $|\xi|\sim \lambda^{-4}|\xi_{1}|^{-3}$
  which yields that $|\xi|\sim \lambda^{-4}|\xi_{1}|^{-3}\sim \lambda^{-4}2^{-3j_{1}}.$
Cases $2^{k_{1}}= 2^{k_{\rm max}}$ and $2^{k_{2}}= 2^{k_{\rm max}}$
can be proved similarly to case
 $$\left|3\xi\xi_{1}\xi_{2}-\frac{\xi_{1}^{2}+\xi_{1}\xi_{2}+\xi_{2}^{2}}{\lambda^{4}\xi\xi_{1}\xi_{2}}
\right|<\frac{|\xi\xi_{1}\xi_{2}|}{4}$$ of (vi) of Lemma 3.1 in this paper.

 \noindent When the subcase (c)  $2^{k}=2^{k_{max}}$ is valid,  we have that $2^{k}\geq C2^{2j_{1}}|\xi|$ and
  $|\xi|\left\langle\sigma^{\lambda}
\right\rangle^{-1}\leq C2^{-2j_{1}}$ or $|\xi|\leq C2^{k-2j_{1}}$ in $B_{k}.$

\noindent In this subcase, we consider two situations
\begin{eqnarray*}
2^{k}\sim {\rm max}\left\{2^{k_{1}},2^{k_{2}}\right\},
 2^{k}\gg{\rm max}\left\{2^{k_{1}},2^{k_{2}}\right\},
\end{eqnarray*}
respectively.

\noindent When $ 2^{k}\sim
{\rm max}\left\{2^{k_{1}},2^{k_{2}}\right\},$
this situation can be proved similarly to case $2^{k_{1}}\sim 2^{k_{\rm max}}$
 or $2^{k_{2}}\sim 2^{k_{\rm max}}$.

\noindent Now we consider the situation  $2^{k}\gg{\rm max}\left\{2^{k_{1}},2^{k_{2}}\right\}$.

\noindent By using the H\"older inequality with respect to $\xi$,
 the Young inequality with respect to $\xi,\tau$  and (\ref{2.013}),
we infer  that
\begin{eqnarray}
&&\left\||\xi|\langle\xi\rangle^{-\frac{3}{4}}\left\langle\sigma^{\lambda}
\right\rangle^{-1}f*g\right\|_{L_{\xi}^{2}
L_{\tau}^{1}(A_{0})}\leq C2^{-2j_{1}}\|f*g\|_{L_{\xi}^{\infty}L_{\tau}^{1}
(A_{0})}\nonumber\\&&\leq C2^{-\frac{j_{1}}{2}}
\left\|\langle\xi\rangle^{-\frac{3}{4}}f\right\|_{L_{\xi}^{2}L_{\tau}^{1}}
\left\|\langle\xi\rangle^{-\frac{3}{4}}g\right\|_{L_{\xi}^{2}L_{\tau}^{1}}
\leq C2^{-\frac{j_{1}}{2}}\|f\|_{\hat{X_{\lambda}}}
\|g\|_{\hat{X_{\lambda}}}\label{3.021}.
\end{eqnarray}
When $(\tau,\xi)\in D^{\prime}$,  note that $\xi_{1}\xi_{2}<0$ in this case, we conclude that
 $K_{1}\sim 2^{j_{1}}$.

\noindent If (\ref{2.02}) is valid, then Lemma 2.1 is valid.
\noindent By  using the fact that $|\xi|^{\frac{1}{2}}\left\langle\sigma^{\lambda}
\right\rangle^{-\frac{1}{2}}\leq C2^{-j_{1}}$   and   Lemma 2.1, we have
\begin{eqnarray}
&&\left\|\xi\left\langle\sigma^{\lambda}\right\rangle^{-\frac{1}{2}}
\langle\xi\rangle^{-3/4}f*g\right\|_{L^{2}_{\xi\tau}}
\leq CK_{1}^{-\frac{1}{2}}2^{-j_{1}}\|f\|_{\hat{X_{\lambda}}^{0,\frac{1}{2},1}}
\|g\|_{\hat{X_{\lambda}}^{0,\frac{1}{2},1}}\nonumber\\
&&\leq C2^{-\frac{3j_{1}}{2}}\|f\|_{\hat{X_{\lambda}}^{0,\frac{1}{2},1}}
\|g\|_{\hat{X_{\lambda}}^{0,\frac{1}{2},1}}\leq
C\|f\|_{\hat{X_{\lambda}}^{-\frac{3}{4},\frac{1}{2},1}}
\|g\|_{\hat{X_{\lambda}}^{-\frac{3}{4},\frac{1}{2},1}}\nonumber\\&&
\leq C\|f\|_{\hat{X_{\lambda}}}
\|g\|_{\hat{X_{\lambda}}}.\label{3.022}
\end{eqnarray}
If  (\ref{2.02}) is invalid, we have that $|\xi|\sim\lambda^{-2}|\xi_{1}|^{-1}\sim \lambda^{-2}2^{-j_{1}}$ which yields that
$\left\langle\sigma^{\lambda}\right\rangle^{-\frac{1}{2}}\sim \lambda^{-3} 2^{-\frac{3j_{1}}{2}}.$ By
using the H\"older inequality with respect to $\xi$,  and the Young inequality with respect to $\xi,\tau$,
we obtain  that
\begin{eqnarray}
&&\left\|\xi\left\langle\sigma^{\lambda}\right\rangle^{-\frac{1}{2}}
\langle\xi\rangle^{-\frac{3}{4}}f*g\right\|_{L^{2}_{\xi\tau}}
\leq C2^{-\frac{5j_{1}}{2}}\|f*g\|_{L_{\xi\tau}^{2}}
\leq C2^{-3j_{1}}\|f*g\|_{L_{\xi}^{\infty}L_{\tau}^{2}}
\nonumber\\&&
\leq C2^{-3j_{1}}\|f\|_{L_{\xi}^{2}L_{\tau}^{2}}\|g\|_{L_{\xi}^{2}L_{\tau}^{1}}\leq
C\|f\|_{\hat{X_{\lambda}}^{-\frac{3}{4},\frac{1}{2},1}}
\|g\|_{\hat{X_{\lambda}}^{-\frac{3}{4},\frac{1}{2},1}}\nonumber\\&&
\leq C\|f\|_{\hat{X_{\lambda}}}
\|g\|_{\hat{X_{\lambda}}}.\label{3.023}
\end{eqnarray}
\noindent When $(\tau,\xi)\in D_{1}\cup D_{2}=\left\{(\tau,\xi)\in \R^{2}:|\xi|\leq1,
|\tau|< |\xi|^{-3}\right\},$ we have that
\begin{eqnarray}
|\xi|^{-3}> |\tau|=\left|\sigma^{\lambda}\right|+1-|\xi|^{3}
-|\lambda^{4}\xi|^{-1}-1\geq 2^{k}-|\xi|^{3}-|\lambda^{4}\xi|^{-1}-1.\label{3.024}
\end{eqnarray}
From  (\ref{3.024}), we have that
$
C|\xi|2^{2j_{1}}\leq 2^{k}\leq |\xi|^{-3}+|\xi|^{3}+|\lambda^{4}\xi|^{-1}+1\leq 4|\xi|^{-3}
$
which yields that
$|\xi|\leq C2^{-\frac{j_{1}}{2}}.$  Thus   this case  can be proved
similarly to cases (3.21)-(3.22)
of \cite{LHY}.

\noindent
When $(\tau,\xi)\in D_{3}=\left\{(\tau,\xi)\in \R^{2},
\frac{1}{8}<|\xi|\leq1, |\tau|\geq |\xi|^{-3}\right\}$.
In this case, we see  that $2^{k}\sim \left|3\xi\xi_{1}\xi_{2}
-\frac{\xi_{1}^{2}+\xi_{1}\xi_{2}+\xi_{2}^{2}}
{\lambda^{4}\xi\xi_{1}\xi_{2}}\right|\sim C|\xi|2^{2j_{1}}$
and (\ref{2.02}) is valid, which implies that Lemma 2.1 is true.
Obviously,  $|\xi|^{\frac{1}{4}}2^{-\frac{k}{2}}
\leq C|\xi|^{-\frac{3}{4}}2^{-j_{1}}|\xi|^{\frac{1}{2}}$.
By  using Lemma 2.1, we infer  that
\begin{eqnarray*}
&&\left\|I_{A_{j}}\left\langle  \sigma^{\lambda} \right\rangle^{-1}
\xi f*g\right\|_{\hat{X}_{\lambda}}
\nonumber\\&&\leq C|\xi|^{\frac{1}{4}}\sum\limits_{k=C(In |\xi|+j_{1}ln 4)+O(1))}
2^{-\frac{k}{2}}\|f*g\|_{L^{2}(A_{j}\cap B_{k})}\nonumber\\
&&\leq C|\xi|^{-\frac{3}{4}}2^{-j_{1}}\sum\limits_{k=C(ln|\xi|+j_{1}ln4)+O(1))}
\||\xi|^{\frac{1}{2}}f*g\|_{L^{2}(A_{j}\cap B_{k})}\nonumber\\
&&\leq C|\xi|^{-\frac{3}{4}}2^{-\frac{3}{2}j_{1}}
\|f\|_{\hat{X}_{\lambda}^{0,\frac{1}{2},1}}
\|g\|_{\hat{X}_{\lambda}^{0,\frac{1}{2},1}}\nonumber\\
&&\leq C\|f\|_{\hat{X}_{\lambda}^{-\frac{3}{4},\frac{1}{2},1}}
\|g\|_{\hat{X}_{\lambda}^{-\frac{3}{4},\frac{1}{2},1}}.
\end{eqnarray*}

\noindent(vii) In this case, we consider two subcases
\begin{eqnarray*}
\left|3\xi\xi_{1}\xi_{2}-\frac{\xi_{1}^{2}+\xi_{1}\xi_{2}
+\xi_{2}^{2}}{\lambda^{4}\xi\xi_{1}\xi_{2}}
\right|<\frac{|\xi\xi_{1}\xi_{2}|}{4}, \;\;\;  \left|3\xi\xi_{1}\xi_{2}-\frac{\xi_{1}^{2}
+\xi_{1}\xi_{2}+\xi_{2}^{2}}{\lambda^{4}\xi\xi_{1}\xi_{2}}
\right|\geq\frac{|\xi\xi_{1}\xi_{2}|}{4}.
\end{eqnarray*}
(vii-1) When $\left|3\xi\xi_{1}\xi_{2}-\frac{\xi_{1}^{2}+\xi_{1}\xi_{2}
+\xi_{2}^{2}}{\lambda^{4}\xi\xi_{1}\xi_{2}}
\right|<\frac{|\xi\xi_{1}\xi_{2}|}{4}$, we have
$|\xi_{2}|\sim  \lambda^{-2}|\xi_{1}|^{-1}\sim\lambda^{-2}2^{-j_{1}}.$

\noindent When $(\tau_{2},\xi_{2})\in D^{\prime}\cup D_{3},$ we have
that $|\tau_{2}|\geq |\xi_{2}|^{-3}
\geq C\lambda^{6}2^{3j_{1}},$ which yields that
\begin{eqnarray}
\left\langle\sigma_{2}^{\lambda}\right
\rangle^{-\frac{1}{2}}\sim C\lambda ^{-3}2^{-\frac{3j_{1}}{2}}.\label{3.025}
\end{eqnarray}
 By using the H\"older inequality  with respect to  $\xi$,   the Young
 inequality with respect to $\xi,\tau$,  together with  (\ref{2.013}) and  (\ref{3.025}),  we have
 \begin{eqnarray}
&&\left\|\chi_{A_{j}}\langle\xi\rangle^{-\frac{3}{4}}\left\langle\sigma^{\lambda}\right\rangle
^{-1}\xi f*g\right\|_{\hat{X_{\lambda}}^{0,\frac{1}{2},1}}\leq
C2^{\frac{j_{1}}{4}}\sum_{k\geq0}2^{-\frac{k}{2}}\|f*g\|_{L_{\xi\tau}^{2}}\nonumber\\&&
\leq C2^{-\frac{3j_{1}}{4}}\|f*\left(\left\langle\sigma^{\lambda}\right
\rangle^{\frac{1}{2}}g\right)\|_{L_{\xi}^{\infty}L_{\tau}^{2}}
\leq C2^{-\frac{3j_{1}}{4}}\|f\|_{L_{\xi}^{2}L_{\tau}^{1}}
\|g\|_{\hat{X_{\lambda}}^{-\frac{3}{4},\frac{1}{2}}}\nonumber\\
&&\leq C\|f\|_{\hat{X_{\lambda}}^{-\frac{3}{4},\frac{1}{2},1}}
\|g\|_{\hat{X_{\lambda}}^{-\frac{3}{4},\frac{1}{2}}}\leq C\|f\|_{\hat{X_{\lambda}}}
\|g\|_{\hat{X_{\lambda}}}.\label{3.026}
\end{eqnarray}
\noindent When $(\tau_{2},\xi_{2})\in  D_{1}\cup D_{2}$,
(\ref{2.01}) or (\ref{2.02}) is valid. Thus,  Lemma 2.1 is valid.
By using Lemma 2.1, we have
\begin{eqnarray}
&&\left\|I_{A_{j}}\langle\xi\rangle^{-\frac{3}{4}}\left\langle\sigma^{\lambda}\right
\rangle^{-1}\xi f*g\right\|_{\hat{X_{\lambda}}^{0,\frac{1}{2},1}}\nonumber\\&&\leq C2^{-\frac{j_{1}}{4}}
\sum_{k\geq 0}2^{-\frac{k}{2}}\||\xi|^{\frac{1}{2}}f*g\|_{L_{\xi\tau}^{2}}
\leq C2^{-\frac{3j_{1}}{4}}
\|f\|_{\hat{X_{\lambda}}^{0,\frac{1}{2},1}}\|g\|_{\hat{X_{\lambda}}^{0,\frac{1}{2},1}}\nonumber\\
&&\leq C\|f\|_{\hat{X_{\lambda}}^{-\frac{3}{4},\frac{1}{2},1}}
\|g\|_{\hat{X_{\lambda}}^{-\frac{3}{4},\frac{1}{2}}}\leq C\|f\|_{\hat{X_{\lambda}}}
\|g\|_{\hat{X_{\lambda}}}.\label{3.027}
\end{eqnarray}
By using (\ref{2.013}), combining (\ref{3.026})  with  (\ref{3.027}),  we have
\begin{eqnarray}
&&\left\|\chi_{A_{j}}\langle\xi\rangle^{-\frac{3}{4}}\left\langle\sigma^{\lambda}\right
\rangle^{-1}\xi f*g\right\|_{L_{\xi}^{2}L_{\tau}^{1}}\nonumber\\&&
\leq C \left\|I_{A_{j}}\langle\xi\rangle^{-\frac{3}{4}}\left\langle\sigma^{\lambda}\right
\rangle^{-1}\xi f*g\right\|_{\hat{X_{\lambda}}^{0,\frac{1}{2},1}}\nonumber\\&&
\leq C\|f\|_{\hat{X_{\lambda}}}
\|g\|_{\hat{X_{\lambda}}}.\label{3.028}
\end{eqnarray}
\noindent (vii-2) When $\left|3\xi\xi_{1}\xi_{2}-\frac{\xi_{1}^{2}+\xi_{1}\xi_{2}
+\xi_{2}^{2}}{\lambda^{4}\xi\xi_{1}\xi_{2}}
\right|\geq\frac{|\xi\xi_{1}\xi_{2}|}{4}$,
we have    $2^{k_{\rm max}}:=2^{{\rm max}\{k,k_{1},k_{2}\}}\geq C|\xi_{2}|2^{2j_{1}}.$
In this subcase, we consider two situations
\begin{eqnarray*}
&&{\rm (a)}:\left|1+\frac{4}{3\lambda ^{4}\xi^{2}\xi_{1}\xi_{2}}\right|
            > \frac{1}{2},\\
&&{\rm (b)}:\left|1+\frac{4}{3\lambda ^{4}\xi^{2}\xi_{1}\xi_{2}}\right|
            \leq\frac{1}{2},
\end{eqnarray*}
respectively.

\noindent (a): When $\left|1+\frac{4}{3\lambda ^{4}\xi^{2}\xi_{1}\xi_{2}}\right|
           \leq \frac{1}{2},$ we  have
\begin{eqnarray}
 |\xi_{2}|\sim \lambda^{-4}|\xi_{1}|^{-3}\sim\lambda^{-4}2^{-3j_{1}}.\label{3.029}
\end{eqnarray}
\noindent By using the  Young inequality with respect to $\xi,\tau$,  and  (\ref{2.013}) and (\ref{3.029}), we have
\begin{eqnarray}
&&\left\|\chi_{A_{j}}\langle\xi\rangle^{-\frac{3}{4}}\left\langle\sigma^{\lambda}\right
\rangle^{-1}\xi f*g\right\|_{L_{\xi}^{2}L_{\tau}^{1}}\leq
C2^{\frac{j_{1}}{4}}\sum_{k\geq 0}2^{-\frac{k}{2}}\|f*g\|_{L_{\xi\tau}^{2}}\nonumber\\&&
\leq C2^{\frac{j_{1}}{4}}\|f*g\|_{L_{\xi}^{2}L_{\tau}^{2}}
\leq C2^{\frac{j_{1}}{4}}\|f\|_{L_{\xi}^{2}L_{\tau}^{1}}\|g\|_{L_{\xi}^{1}L_{\tau}^{2}}
\nonumber\\&&\leq C2^{-\frac{5j_{1}}{4}}\|f\|_{L_{\xi}^{2}L_{\tau}^{1}}\|g\|_{L_{\xi}^{2}L_{\tau}^{2}}
\nonumber\\
&&\leq C\|f\|_{\hat{X_{\lambda}}^{-\frac{3}{4},\frac{1}{2},1}}
\|g\|_{\hat{X_{\lambda}}^{-\frac{3}{4},\frac{1}{2}}}\leq C\|f\|_{\hat{X_{\lambda}}}
\|g\|_{\hat{X_{\lambda}}}\label{3.030}
\end{eqnarray}
    and
\begin{eqnarray}
&&\left\|I_{A_{j}}\langle\xi\rangle^{-\frac{3}{4}}\left\langle\sigma^{\lambda}\right
\rangle^{-1}\xi f*g\right\|_{\hat{X_{\lambda}}^{0,\frac{1}{2},1}}\nonumber\\&&\leq
C2^{\frac{j_{1}}{4}}\sum_{k\geq 0}2^{-\frac{k}{2}}\|f*g\|_{L_{\xi\tau}^{2}}\nonumber\\&&
\leq C2^{\frac{j_{1}}{4}}\|f*g\|_{L_{\xi}^{2}L_{\tau}^{2}}
\leq C2^{\frac{j_{1}}{4}}\|f\|_{L_{\xi}^{2}L_{\tau}^{1}}\|g\|_{L_{\xi}^{1}L_{\tau}^{2}}
\nonumber\\&&\leq C2^{\frac{-5j_{1}}{4}}\|f\|_{L_{\xi}^{2}L_{\tau}^{1}}\|g\|_{L_{\xi}^{2}L_{\tau}^{2}}
\nonumber\\
&&\leq C\|f\|_{\hat{X_{\lambda}}^{-\frac{3}{4},\frac{1}{2},1}}
\|g\|_{\hat{X_{\lambda}}^{-\frac{3}{4},\frac{1}{2}}}\leq C\|f\|_{\hat{X_{\lambda}}}
\|g\|_{\hat{X_{\lambda}}}.\label{3.031}
\end{eqnarray}
\noindent (b): When $\left|1+\frac{4}{3\lambda ^{4}\xi^{2}\xi_{1}\xi_{2}}\right|
           >\frac{1}{2},$ we consider three scenarios
\begin{eqnarray*}
(b1):2^{k}= 2^{{\rm max}};\;\;   (b2):2^{k_{1}}= 2^{{\rm max}};\;\;    (b3):2^{k_{2}}= 2^{{\rm max}},
\end{eqnarray*}
 respectively.

\noindent (b1)  $2^{k}= 2^{{\rm max}}$ can be proved similarly to
 case $2^{k}= 2^{{\rm max}}$   of (v) in \cite{Kis} since Lemma 2.1 is valid in this case.

\noindent (b2) $2^{k_{1}}= 2^{{\rm max}}$ can be proved similarly
to case $2^{k_{1}}= 2^{{\rm max}}$   of (v) in \cite{Kis} with the aid of (\ref{2.010}).

\noindent (b3) $2^{k_{2}}= 2^{{\rm max}}$. In this scenario,
we consider
\begin{eqnarray}
2^{k_{2}}\gg{\rm max}\left\{2^{k},2^{k_{1}}\right\},2^{k_{2}}\sim {\rm max}\left\{2^{k},2^{k_{1}}\right\},
\end{eqnarray}
respectively.

\noindent When $2^{k_{2}}\sim{\rm max}\left\{2^{k},2^{k_{1}}\right\}$,
 this case can be proved similarly
to case $2^{k}= 2^{{\rm max}}$ or case
$2^{k_{1}}= 2^{{\rm max}}$  of (v) in \cite{Kis}.

\noindent When $2^{k_{2}}\gg{\rm max}\left\{2^{k},2^{k_{1}}\right\}$,
we have that
\begin{eqnarray}
2^{k_{2}}\sim\left|3\xi\xi_{1}\xi_{2}-\frac{\xi_{1}^{2}
+\xi_{1}\xi_{2}+\xi_{2}^{2}}{\lambda^{4}\xi\xi_{1}\xi_{2}}
\right|.\label{3.032}
\end{eqnarray}

\noindent We consider $(\tau_{2},\xi_{2})\in D^{\prime}$
 and $(\tau_{2},\xi_{2})\in D_{1}\cup D_{2}\cup D_{3},$ respectively.

\noindent
When $(\tau_{2},\xi_{2})\in D^{\prime}$, from (\ref{3.025})-(\ref{3.028})
of \cite{LHY}, we know that
\begin{eqnarray}
|\xi\xi_{1}\xi_{2}|\geq\frac{\xi_{1}^{2}+\xi_{1}\xi_{2}
+\xi_{2}^{2}}{\lambda^{4}|\xi\xi_{1}\xi_{2}|}.\label{3.033}
\end{eqnarray}
Combining (\ref{3.032}) with (\ref{3.033}), we know that
\begin{eqnarray}
\frac{1}{2}|\xi_{2}|^{-3}\leq  2^{k_{2}}\sim\left|3\xi\xi_{1}\xi_{2}
-\frac{\xi_{1}^{2}+\xi_{1}\xi_{2}+\xi_{2}^{2}}{\lambda^{4}\xi\xi_{1}\xi_{2}}
\right|\sim |\xi\xi_{1}\xi_{2}|\sim |\xi_{2}|2^{2j_{1}}.\label{3.034}
\end{eqnarray}
From (\ref{3.034}), we have that
\begin{eqnarray}
|\xi_{2}|\geq C2^{-\frac{j_{1}}{2}}.\label{3.035}
\end{eqnarray}
Thus, we have that $C2^{\frac{3j}{2}}\leq 2^{k_{2}}\leq C2^{2j}$.
 In this case, we have that
 \begin{eqnarray}
\left|1-\frac{4}{3\lambda ^{4}\xi\xi_{1}\xi_{2}^{2}}\right|
            > \frac{1}{2}.\label{3.036}
\end{eqnarray}
\noindent Since  $\left|1+\frac{4}{3\lambda ^{4}\xi^{2}\xi_{1}\xi_{2}}\right|
           >\frac{1}{2}$ and (\ref{3.036}) is valid,  we know that Lemmas 2.1, 2.2 are valid.
This case can be proved similarly to case $(\tau_{2},\xi_{2})\in D$ of (v) of
 Proposition 3.4 in \cite{Kis}.

\noindent When $(\tau_{2},\xi_{2})\in D_{1}\cup D_{2}\cup D_{3},$  by using
 $|\xi_{2}|\leq C2^{k_{2}}2^{-2j}$  together with    (\ref{2.010}),
  (\ref{2.012}) and  (\ref{2.011}),
we obtain that
\begin{eqnarray}
&&\left\|\chi_{A_{j}}\langle\xi\rangle^{-\frac{3}{4}}\left\langle\sigma^{\lambda}\right
\rangle^{-1}\xi f*g\right\|_{L_{\xi}^{2}L_{\tau}^{1}}\leq
C2^{j_{1}}\sum_{k\geq 0}2^{-\frac{k}{2}}\|(\langle \xi\rangle ^{-\frac{3}{4}}f)*g\|_{L_{\xi\tau}^{2}}\nonumber\\&&
\leq C2^{j_{1}}\|(\langle \xi\rangle ^{-\frac{3}{4}}f)*g\|_{\hat{X_{\lambda}}^{0,-\frac{1}{2}+\epsilon}}\nonumber\\&&
\leq C2^{j_{1}}\|(\langle \xi\rangle ^{-\frac{3}{4}}f)*g\|_{L_{\xi}^{2}L_{\tau}^{4}}\nonumber\\&&
\leq C2^{j}\|\langle \xi\rangle ^{-\frac{3}{4}}f\|_{L_{\xi}^{2}L_{\tau}^{\frac{4}{3}}}
\sum\limits_{k_{2}}\|g\|_{L_{\xi}^{1}L_{\tau}^{2}
(B_{k_{2}}\cap \{|\xi|\leq C2^{k_{2}-2j}\})}\nonumber\\&&\leq C\|f\|_{\hat{X_{\lambda}}^{-\frac{3}{4},\frac{1}{2}}}
\|g\|_{\hat{X_{\lambda}}^{-\frac{3}{4},\frac{1}{2},1}}\nonumber\\&&\leq C\|f\|_{\hat{X_{\lambda}}}
\|g\|_{\hat{X_{\lambda}}}.\label{3.037}
\end{eqnarray}

\noindent(viii) This case can be proved similarly to case (vii) due to the symmetry.

This completes   the proof of Lemma 3.1.

\noindent{\bf Remark 7.}  Cases (vi) and (vii) are the most difficult parts   in
  proving Lemma 3.1.  The proof   of   Lemma 3.1 demonstrates
that the structure of the Ostrovsky equation with positive dispersion are much more complicated
than that   with negative dispersion and of   the KdV equation. Thus, we cannot completely follow the method of \cite{Kis,LHY}.  It is necessary that we should pay more attention to the   structure of the
Ostrovsky equation with positive dispersion, as the phase function $\phi^{\lambda}(\xi)$ and the resonant function $\left|\sigma^{\lambda}-\sigma_{1}^{\lambda}
-\sigma_{2}^{\lambda}\right|
=\left|3\xi\xi_{1}\xi_{2}-\frac{\xi_{1}^{2}+\xi_{1}\xi_{2}
+\xi_{2}^{2}}{\lambda ^{4}\xi\xi_{1}\xi_{2}}\right|$ of the
Ostrovsky equation with positive dispersion are much more complicated than those    with negative dispersion and of  KdV equation.

\begin{Lemma}\label{Lem3.2}
Let  $u,v\in X_{\lambda}$.  Then
\begin{eqnarray}
&&\hspace{-1.5cm}\left\|\mathscr{F}^{-1}\left[\left\langle\sigma^{\lambda}
\right\rangle^{-1}\mathscr{F}\left[\partial_{x}(uv)\right]\right]
\right\|_{X_{\lambda}}+
\left\|\mathscr{F}^{-1}\left[\left\langle\sigma^{\lambda}\right\rangle^{-1}
\mathscr{F}\left[\partial_{x}(uv)\right]\right]
\right\|_{Y}\leq C\|u\|_{X_{\lambda}}\|v\|_{X_{\lambda}}.\label{3.038}
\end{eqnarray}
\end{Lemma}
\noindent{\bf Proof.}
Combining the technique of Lemma 3.2 of \cite{LHY} with Lemma 3.1,  we have Lemma 3.2.

This completes the proof of Lemma 3.2.

\bigskip
\bigskip

\noindent {\large\bf 4. Proof of Theorem  1.1: Well-posedness for $s=-\frac34$}

\setcounter{equation}{0}

 \setcounter{Theorem}{0}

\setcounter{Lemma}{0}

\setcounter{section}{4}

In this section, we use  Lemmas 2.4, 2.5, 3.2 to  prove Theorem 1.1,
 the local well-posedness   for $s=-\frac34$.  Without loss of generality,  we assume that $T=1$.

 The  solution to the  Cauchy problem for (\ref{1.03}) can be formally rewritten as follows:
\begin{eqnarray}
u(t)=e^{-t(-\partial_{x}^{3}-\partial_{x}^{-1})}u(x,0)+\frac{1}{2}\int_{0}^{t}e^{-(t-s)
(-\partial_{x}^{3}-\partial_{x}^{-1})}\partial_{x}(u^{2})ds.\label{4.01}
\end{eqnarray}
We define
\begin{eqnarray}
\Phi(u):=\psi(t)e^{-t(-\partial_{x}^{3}-\partial_{x}^{-1})}u(x,0)+\frac{1}{2}\psi(t)\int_{0}^{t}e^{-(t-s)
(-\partial_{x}^{3}-\partial_{x}^{-1})}\partial_{x}(u^{2})ds.\label{4.02}
\end{eqnarray}
By taking advantage of  Lemmas 2.4, 2.5, 3.2  with $\lambda =1$, we derive that
\begin{eqnarray}
\|\Phi (u)\|_{X_{1},1}+\sup\limits_{-1\leq t\leq1}\|\Phi(u)\|_{H^{-\frac34}(\SR)}
\leq C\|u(x,0)\|_{H^{-\frac34}(\SR)}+C\|u\|_{X_{1},1}^{2}.  \label{4.03}
\end{eqnarray}
If $\|u(x,0)\|_{H^{-\frac34}}$ is sufficiently small,  then  $\Phi(u)$ is
a contraction mapping on  the closed ball
\begin{eqnarray}
B=\left\{u\in X_{_{1},1}: \|u\|_{X_{1},1}\leq 2C\|u_{0}\|_{H^{-\frac34}(\SR)}\right\}.\label{4.04}
\end{eqnarray}
Thus $\Phi$ has a fixed point $u$ on closed ball $B$  and
consequently, the Cauchy problem for (\ref{1.03}) possesses a local solution on $[-1,1]$.
For the uniqueness of the solutions, we refer the readers to   \cite{MT,Kis}.

\noindent
Now we prove that  $v\in C([-1,1]; H^{-\frac{3}{4}}(\R))$, by using Lemmas 2.4, 3.2, for $t\in [-1,1]$, we can easily obtain
\begin{eqnarray}
&&\|u(\cdot,t)\|_{H^{-\frac{3}{4}}(\SR)}\nonumber\\&&\leq C\left\|e^{-t(-\partial_{x}^{3}-\partial_{x}^{-1})}
u(x,0)\right\|_{H^{-\frac{3}{4}}(\SR)}+C\left\|\int_{0}^{t}e^{-(t-s)
(-\partial_{x}^{3}-\partial_{x}^{-1})}\partial_{x}(u^{2})ds\right\|_{H^{-\frac{3}{4}}(\SR)}\nonumber\\&&
\leq C\|u(\cdot,0)\|_{H^{-\frac{3}{4}}(\SR)}+C\left\|\langle \xi\rangle^{-\frac{3}{4}} \left\langle \sigma^{1}\right\rangle^{-1}\xi\mathscr{F}(u^{2})\right\|_{L_{\xi}^{2}L_{\tau}^{1}}\nonumber\\
&&\leq C\left[\|u(\cdot,0)\|_{H^{-\frac{3}{4}}(\SR)}+\|u\|_{X_{1},1}^{2}\right],\label{4.05}
\end{eqnarray}
then, we have that  $u\in L^{\infty}([-1,1]; H^{-\frac{3}{4}}(\R)).$
Let $u(x,t_{1})$ and $u(x,t_{2})$ be the solutions to (\ref{1.03})-(\ref{1.02}) at the moment $t_{1}$, $t_{2}$, respectively.
By using  a proof similar to (\ref{2.019})-(\ref{2.021}), we know that for $\forall \epsilon_{1}>0,$
there exist sufficiently large constant $A_{1}>0$ and $\delta_{1}>0$,  which satisfies
 $0<\delta_{1}<\frac{\epsilon_{1}}{2A_{1}\|u\|_{X_{\lambda}^{-\frac34,\frac{1}{2},1}}},$  and when
$|t_{1}-t_{2}|<\delta_{1},$ we have
\begin{eqnarray}
\|u(\cdot,t_{1})-u(\cdot,t_{2})\|_{H^{-\frac34}(\SR)}<\epsilon_{1}.\label{4.06}
\end{eqnarray}
Combining the conclusion  that $u\in L^{\infty}([-1,1]; H^{-\frac{3}{4}}(\R))$
 with (\ref{4.07}), since $\epsilon_1,\epsilon_2$ are arbitrarily small,
  we have that
$u\in C([-1,1]; H^{-\frac{3}{4}}(\R))$.
Now we prove that  the Lipschitz dependence of solutions
on the initial data holds.
We assume that $u_{1}(x,t),u_{2}(x,t)$  are solutions to (\ref{1.03}), (\ref{1.02})
corresponding to the initial data $u_{1}(x,0),u_{2}(x,0)$, respectively.
Then, by exploiting Lemmas 2.4, 2.5, 3.2 and (\ref{4.04}), since $\|u_{j}(x,0)\|_{H^{-\frac34}}$
 is sufficiently small with $j=1,2,$ we derive that
\begin{eqnarray}
&&\|u_{1}-u_{2}\|_{X_{1},1}+\sup\limits_{-1\leq t\leq1}\|u_{1}-u_{2}\|_{H^{-\frac34}(\SR)}\nonumber\\&&
\leq C\|u_{1}(\cdot,0)-u_{2}(\cdot,0)\|_{H^{-\frac34}(\SR)}+C\|u_{1}-u_{2}\|_{X_{1},1}
\left[\|u_{1}\|_{X_{1},1}^{2}+\|u_{2}\|_{X_{1},1}\right]\nonumber\\
&&\leq C\|u_{1}(\cdot,0)-u_{2}(\cdot,0)\|_{H^{-\frac34}(\SR)}+C\|u_{1}-u_{2}\|_{X_{1},1}
\left[\|u_{1}(\cdot,0)\|_{H^{-\frac{3}{4}}(\SR)}+\|u_{2}(\cdot,0)\|_{H^{-\frac{3}{4}}(\SR)}\right]\nonumber\\
&&\leq C\|u_{1}(\cdot,0)-u_{2}(\cdot,0)\|_{H^{-\frac34}(\SR)}+\frac{1}{2}\|u_{1}-u_{2}\|_{X_{1},1}.  \label{4.07}
\end{eqnarray}
From (\ref{4.07}), we have that
\begin{eqnarray}
\frac{1}{2}\|u_{1}-u_{2}\|_{X_{1},1}+\sup\limits_{-1\leq t\leq1}\|u_{1}(\cdot,t)-u_{2}(\cdot,t)\|_{H^{-\frac34}(\SR)}
\leq C\|u_{1}(\cdot,0)-u_{2}(\cdot,0)\|_{H^{-\frac34}(\SR)}.  \label{4.08}
\end{eqnarray}
From (\ref{4.08}), we have that
\begin{eqnarray}
\|u_{1}-u_{2}\|_{X_{1},1}+2\sup\limits_{-1\leq t\leq1}\|u_{1}(\cdot,t)-u_{2}(\cdot,t)\|_{H^{-\frac34}(\SR)}
\leq 2C\|u_{1}(\cdot,0)-u_{2}(\cdot,0)\|_{H^{-\frac34}(\SR)}.  \label{4.09}
\end{eqnarray}
Thus, we establish the conclusion that the Lipschitz dependence of solutions
on the initial data holds.

\noindent
When  $\|u(\cdot,0)\|_{H^{-\frac34}(\SR)}$ is arbitrarily large, it is easily checked that $u^{\lambda}(x,t)=\lambda^{-2}u\left(\frac{x}{\lambda},\frac{t}{\lambda^{3}}\right)$ is
 the solution to (\ref{1.04})-(\ref{1.05}).
By using a direct computation, we have that
 \begin{eqnarray*}
 \|u^{\lambda}(\cdot,0)\|_{H^{-\frac34}(\SR)}\leq C\lambda^{-\frac34}\|u(x,0)\|_{H^{-\frac34}(\SR)}.
\end{eqnarray*}
Taking $\lambda$ sufficiently large, then we have that
$\|u^{\lambda}(\cdot,0)\|_{H^{-\frac34}(\SR)}$  is   sufficiently small, this case can be proved similarly
to case that $\|u(\cdot,0)\|_{H^{-\frac34}(\SR)}$  is   sufficiently small.

This completes the proof of Theorem 1.1.

\bigskip
\bigskip

\noindent{\large\bf 5. Bilinear estimate for $s>-\frac34$
 and the local well-posedness for $s>-\frac34$}

\setcounter{equation}{0}

 \setcounter{Theorem}{0}

\setcounter{Lemma}{0}

 \setcounter{section}{5}

We now prove Theorem 1.2.  More precisely, we use  Lemmas 2.7, 2.9  to
establish the bilinear estimate  for  $s > -\frac34$.
It suffices to prove the following Lemma 5.1 for arbitrarily  positive $\epsilon$.

\begin{Lemma}\label{Lem5.1}
Let  $u_{j}\in X_{\lambda}^{s,\frac{1}{2}+\epsilon}$
 with $s\geq-\frac{3}{4}+6\epsilon$
 and $j=1,2.$ Then, we have
\begin{eqnarray}
\left\|\partial_{x}\left(\prod_{j=1}^{2}u_{j}\right)
\right\|_{X_{\lambda}^{s,-\frac{1}{2}+2\epsilon}}\leq
 C\prod_{j=1}^{2}\|u_{j}\|
 _{X_{\lambda}^{s,\frac{1}{2}+\epsilon}}.\label{5.01}
\end{eqnarray}
\end{Lemma}
\noindent {\bf Proof.} To prove (\ref{5.01}), it suffices by duality  to show that
\begin{eqnarray}
\int_{\SR^{2}}\bar{u}(x,t)\partial_{x}\left(\prod_{j=1}^{2}u_{j}\right)dxdt\leq
 C\left[\prod_{j=1}^{2}\|u_{j}\|_{X_{\lambda}^{s,\frac{1}{2}+\epsilon}}\right]
 \|u\|_{X_{\lambda}^{-s,\frac{1}{2}-2\epsilon}}.\label{5.02}
\end{eqnarray}
Let
\begin{eqnarray}
&&\hspace{-1cm}F(\xi,\tau)=\langle\xi\rangle^{-s}\langle \sigma^{\lambda}\rangle^{\frac{1}{2}-2\epsilon}\mathscr{F}u(\xi,\tau),
F_{j}(\xi_{j},\tau_{j})=\langle\xi_{j}\rangle^{s}\langle \sigma_{j}^{\lambda}\rangle^{\frac{1}{2}+\epsilon}
\mathscr{F}u_{j}(\xi_{j},\tau_{j})(j=1,2).\label{5.03}
\end{eqnarray}
To obtain (\ref{5.02}),  as seen in  (\ref{5.03}),  it suffices to prove that
\begin{eqnarray}
&&\hspace{-1.2cm}\int_{\SR^2}\int_{\tiny\begin{array}{cc}\xi=\sum\limits_{j=1}^{2}\xi_j\\
\tau =\sum\limits_{j=1}^{2}\tau_j\end{array}}\frac{|\xi|\langle\xi\rangle^{s}
F(\xi,\tau)\prod\limits_{j=1}^{2}F_{j}(\xi_{j},\tau_{j})}{\langle\sigma_{j}^{+}\rangle^{\frac{1}{2}-2\epsilon}
\prod\limits_{j=1}^{2}\langle\xi_{j}\rangle^{s}\langle\sigma_{j}^{+}
\rangle^{\frac{1}{2}+\epsilon}}d\xi_{1}d\tau_{1}d\xi d\tau\leq C
\|F\|_{L_{\xi\tau}^{2}}\left(\prod_{j=1}^{2}\|F_{j}\|_{L_{\xi\tau}^{2}}\right).\label{5.04}
\end{eqnarray}
Without loss of generality and using the symmetry,  we assume that
$|\xi_{1}|\geq |\xi_{2}|$  and   $F(\xi,\tau)\geq 0,F_j(\xi_{j},\tau_{j})\geq 0(j=1,2).$
We define
\begin{eqnarray*}
\Omega^{*}:=\{(\xi_1,\tau_1,\xi,\tau)\in {\rm R^4},\xi=\sum\limits_{j=1}^{2}\xi_j,\tau=\sum\limits_{j=1}^{2}\tau_j,
 |\xi_2|\leq |\xi_{1}|\}
\end{eqnarray*}
and
\begin{eqnarray*}
&&\hspace{-0.8cm}\Omega_1=\{(\xi_1,\tau_1,\xi,\tau)\in\Omega^{*},
 |\xi_2|\leq |\xi_{1}|\leq 16 \},\\
&&\hspace{-0.8cm} \Omega_2=\{ (\xi_1,\tau_1,\xi,\tau)\in \Omega^{*},
|\xi_1|\geq 16, |\xi_{1}|\geq4|\xi_{2}|\},\\
&&\hspace{-0.8cm}\Omega_3=\{(\xi_1,\tau_1,\xi,\tau)\in \Omega^{*},
|\xi_{1}|\geq 16, |\xi_{2}|\leq |\xi_{1}|\leq 4|\xi_2|,\xi_{1}\xi_{2}<0,|\xi|\leq \frac{|\xi_{2}|}{4}\},\\
&&\hspace{-0.8cm}\Omega_4=\{(\xi_1,\tau_1,\xi,\tau)\in \Omega^{*},
 |\xi_{1}|\geq 16, |\xi_{2}|\leq |\xi_{1}|\leq 4|\xi_2|,\xi_{1}\xi_{2}<0,|\xi|\geq \frac{|\xi_{2}|}{4}\},\\
 &&\hspace{-0.8cm}\Omega_5=\{(\xi_1,\tau_1,\xi,\tau)\in \Omega^{*},
 |\xi_{1}|\geq 16, |\xi_{2}|\leq |\xi_{1}|\leq 4|\xi_2|,\xi_{1}\xi_{2}>0\}.
\end{eqnarray*}
 Note that $\Omega^{*}\subset\bigcup\limits_{j=1}^{5}\Omega_{j}.$
Define
\begin{equation}\label{4.05}
    K_{3}(\xi_{1},\tau_{1},\xi,\tau):=\frac{|\xi|\langle \xi \rangle^{s}}
    {\langle\sigma^{\lambda}\rangle^{\frac{1}{2}-2\epsilon}\prod\limits_{j=1}^{2}\langle \xi_{j}\rangle^{s}
    \langle\sigma_{j}^{\lambda}\rangle^{\frac{1}{2}+\epsilon}}
\end{equation}
and
\begin{eqnarray*}
Int:=\int_{\SR^2}\int_{\tiny\begin{array}{cc}\xi=\sum\limits_{j=1}^{2}\xi_j\\ \tau =\sum\limits_{j=1}^{2}\tau_j\end{array}}K_{3}(\xi_{1},\tau_{1},\xi,\tau)F(\xi,\tau)\prod_{j=1}^{2}F_{j}(\xi_{j},\tau_{j})
d\xi_{1}d\tau_{1}d\xi d\tau.
\end{eqnarray*}
\noindent (1) Region $\Omega_{1}$. In this subregion, we have that
\begin{eqnarray}
 K_{3}(\xi_{1},\tau_{1},\xi,\tau)\leq \frac{C}{\langle\sigma^{\lambda}\rangle^{\frac{1}{2}-2\epsilon}
 \prod\limits_{j=1}^{2}\langle\sigma_{j}^{\lambda}\rangle^{\frac{1}{2}+\epsilon}}.\label{5.06}
\end{eqnarray}
By using (\ref{5.06}),  the Cauchy-Schwartz inequality, Plancherel identity and
 the H\"older inequality as well as (\ref{2.021}),  we obtain  that
\begin{eqnarray*}
&&Int\leq C\int_{\SR^2}\int_{\tiny\begin{array}{cc}\xi=\sum\limits_{j=1}^{2}\xi_j\\
\tau =\sum\limits_{j=1}^{2}\tau_j\end{array}}\frac{
F\prod\limits_{j=1}^{2}F_{j}}{\langle\sigma^{\lambda}\rangle^{\frac{1}{2}-2\epsilon}
\prod\limits_{j=1}^{2}\langle\sigma_{j}^{\lambda}\rangle^{\frac{1}{2}+\epsilon}}
d\xi_{1}d\tau_{1}d\xi d\tau\nonumber\\
&&\leq C \left\|\frac{F}{\langle\sigma^{\lambda}\rangle^{\frac{1}{2}-2\epsilon}}\right\|_{L_{\xi\tau}^{2}}
\left\|\int_{\tiny\begin{array}{cc}\xi=\sum\limits_{j=1}^{2}\xi_j\\ \tau
=\sum\limits_{j=1}^{2}\tau_j\end{array}}\frac{\prod\limits_{j=1}^{2}F_{j}}{\prod\limits_{j=1}^{2}
\langle\sigma_{j}^{\lambda}\rangle^{\frac{1}{2}+\epsilon}}d\xi_{1}d\tau_{1}\right\|_{L_{\xi\tau}^{2}}\nonumber\\
&&\leq C\|F\|_{L_{\xi\tau}^{2}}
\left(\prod\limits_{j=1}^{2}\left\|\mathscr{F}^{-1}\left(\frac{F_{j}}{\langle\sigma_{j}^{\lambda}
\rangle^{\frac{1}{2}+\epsilon}}\right)\right\|_{L_{xt}^{4}}\right)\leq C\|F\|_{L_{\xi\tau}^{2}}
\left(\prod\limits_{j=1}^{2}\|F_{j}\|_{L_{\xi\tau}^{2}}\right).
\end{eqnarray*}
\noindent (2) Region $\Omega_{2}$. If  $|\xi_{2}|\leq 1,$  we infer  that
\begin{eqnarray}
 K_{3}(\xi_{1},\tau_{1},\xi,\tau)\leq \frac{C|\xi|}{\langle\sigma^{\lambda}\rangle^{\frac{1}{2}-2\epsilon}
 \prod\limits_{j=1}^{2}\langle\sigma_{j}^{\lambda}\rangle^{\frac{1}{2}+\epsilon}}\leq \frac{C|\left[\phi^{\lambda}(\xi_{1})\right]^{\prime}-\left[\phi^{\lambda}(\xi_{2})\right]^{\prime}|^{\frac{1}{2}}}
 {\langle\sigma^{\lambda}\rangle^{\frac{1}{2}-2\epsilon}
 \prod\limits_{j=1}^{2}\langle\sigma_{j}^{\lambda}\rangle^{\frac{1}{2}+\epsilon}} .\label{5.07}
\end{eqnarray}
By using (\ref{5.07}),  the Cauchy-Schwartz inequality,  Plancherel identity and
  the H\"older inequality as well as Lemma 2.9,
we conclude  that
\begin{eqnarray*}
&&Int\leq C\int_{\SR^2}\int_{\tiny\begin{array}{cc}\xi=\sum\limits_{j=1}^{2}\xi_j\\
\tau =\sum\limits_{j=1}^{2}\tau_j\end{array}}\frac{|\left[\phi^{\lambda}
(\xi_{1})\right]^{\prime}-\left[\phi^{\lambda}(\xi_{2})\right]^{\prime}|^{\frac{1}{2}}
F\prod\limits_{j=1}^{2}F_{j}}{\langle\sigma^{\lambda}\rangle^{\frac{1}{2}-2\epsilon}
\prod\limits_{j=1}^{2}\langle\sigma_{j}^{\lambda}\rangle^{\frac{1}{2}+\epsilon}}
d\xi_{1}d\tau_{1}d\xi d\tau\nonumber\\
&&\leq C \left\|\mathscr{F}^{-1}\left(\frac{F}
{\langle\sigma^{\lambda}\rangle^{\frac{1}{2}-2\epsilon}}\right)\right\|_{L_{xt}^{2}}
\left\|I^{\frac{1}{2}}\left(\mathscr{F}^{-1}\left(\frac{F_{1}}{\langle \sigma_{1}^{\lambda} \rangle ^{\frac{1}{2}+\epsilon}}\right),\mathscr{F}^{-1}\left(\frac{F_{2}}{\langle \sigma_{2}^{\lambda} \rangle ^{\frac{1}{2}+\epsilon}}\right)\right)\right\|_{L_{xt}^{2}}\nonumber\\
&&\leq C\|F\|_{L_{\xi\tau}^{2}}\left(\prod_{j=1}^{2}\|F_{j}\|_{L_{\xi\tau}^{2}}\right).
\end{eqnarray*}
\noindent If $|\xi_{2}|\geq 1,$ and $s\geq 0$,  then we obtain  that
\begin{eqnarray*}
 K_{3}(\xi_{1},\tau_{1},\xi,\tau)\leq \frac{C|\xi|}{\langle\sigma^{\lambda}\rangle^{\frac{1}{2}-2\epsilon}
 \prod\limits_{j=1}^{2}\langle\sigma_{j}^{\lambda}\rangle^{\frac{1}{2}+\epsilon}}\leq \frac{C|\left[\phi^{\lambda}(\xi_{1})\right]^{\prime}-\left[\phi^{\lambda}(\xi_{2})\right]^{\prime}|^{\frac{1}{2}}}
 {\langle\sigma^{\lambda}\rangle^{\frac{1}{2}-2\epsilon}\prod\limits_{j=1}^{2}
 \langle\sigma_{j}^{\lambda}\rangle^{\frac{1}{2}+\epsilon}} .
\end{eqnarray*}
This case can be proved similarly to case $|\xi_{2}|\leq1$ of $\Omega_{2}$ of Lemma 5.1.

\noindent If $|\xi_{2}|\geq 1$ and $-\frac{3}{4}+6\epsilon\leq s<0$, then
 (\ref{1.017}) is valid.

\noindent  Note that $|\sigma^{\lambda}-\sigma_{1}^{\lambda}-\sigma_{2}^{\lambda}|
=\left|3\xi_{1}\xi_{2}(\xi_{1}+\xi_{2})-\frac{\xi_{1}^{2}+\xi_{1}\xi_{2}+\xi_{2}^{2}}
{\lambda^{4}\xi_{1}\xi_{2}(\xi_{1}+\xi_{2})}\right|\geq C|\xi\xi_{1}\xi_{2}|$.
Thus  one of the following three cases must occur:
\begin{eqnarray}
&&|\sigma^{\lambda}|:={\rm max}\left\{|\sigma^{\lambda}|,
|\sigma_{1}^{\lambda}|,|\sigma_{2}^{\lambda}|\right\}\geq C|\xi\xi_{1}\xi_{2}|,\label{5.08}\\
&&|\sigma_{1}^{\lambda}|:={\rm max}\left\{|\sigma^{\lambda}|,
|\sigma_{1}^{\lambda}|,|\sigma_{2}^{\lambda}|\right\}\geq C|\xi\xi_{1}\xi_{2}|,\label{5.09}\\
&&|\sigma_{2}^{\lambda}|:={\rm max}\left\{|\sigma^{\lambda}|,
|\sigma_{1}^{\lambda}|,|\sigma_{2}^{\lambda}|\right\}\geq C|\xi\xi_{1}\xi_{2}|.\label{5.010}
\end{eqnarray}
When (\ref{5.08})  is valid,  we have that
\begin{eqnarray}
&& K_{3}(\xi_{1},\tau_{1},\xi,\tau)\leq \frac{C|\xi||\xi_{2}|^{-s}}
{\langle\sigma^{\lambda}\rangle^{\frac{1}{2}-2\epsilon}
\prod\limits_{j=1}^{2}\langle\sigma_{j}^{\lambda}\rangle^{\frac{1}{2}+\epsilon}}\leq
  C\frac{|\xi|^{2\epsilon}|\xi_{2}|^{\frac{1}{4}-4\epsilon}}{\prod\limits_{j=1}^{2}
 \langle\sigma_{j}^{\lambda}\rangle^{\frac{1}{2}+\epsilon}}\nonumber\\&&\leq C
 \frac{|\left[\phi^{\lambda}(\xi_{1})\right]^{\prime}-\left[\phi^{\lambda}(\xi_{2})\right]^{\prime}|^{\frac{1}{2}}}
 {\prod\limits_{j=1}^{2}\langle\sigma_{j}^{\lambda}\rangle^{\frac{1}{2}+\epsilon}} .\label{5.011}
\end{eqnarray}
Consequently, by using the Plancherel identity, the Cauchy-Schwartz inequality, the H\"older inequality and Lemma 2.9,
we infer    from (\ref{5.011})  that
\begin{eqnarray}
&&Int\leq C\int_{\SR^2}\int_{\tiny\begin{array}{cc}\xi=\sum\limits_{j=1}^{2}\xi_j\\ \tau =\sum\limits_{j=1}^{2}\tau_j\end{array}}\frac{|\left[\phi^{\lambda}(\xi_{1})
\right]^{\prime}-\left[\phi^{\lambda}(\xi_{2})\right]^{\prime}|^{\frac{1}{2}}
F\prod\limits_{j=1}^{2}F_{j}}{\prod\limits_{j=1}^{2}\langle\sigma_{j}^{\lambda}\rangle^{\frac{1}{2}+\epsilon}}
d\xi_{1}d\tau_{1}d\xi d\tau\nonumber\\
&&\leq C \|F\|_{L_{\xi\tau}^{2}}
\left\|I^{\frac{1}{2}}\left(\mathscr{F}^{-1}\left(\frac{F_{1}}{\langle \sigma_{1}^{\lambda} \rangle ^{\frac{1}{2}+\epsilon}}\right),\mathscr{F}^{-1}\left(\frac{F_{2}}{\langle \sigma_{2}^{\lambda} \rangle ^{\frac{1}{2}+\epsilon}}\right)\right)\right\|_{L_{xt}^{2}}\nonumber\\
&&\leq C\|F\|_{L_{\xi\tau}^{2}}\left(\prod_{j=1}^{2}\|F_{j}\|_{L_{\xi\tau}^{2}}\right).\label{5.012}
\end{eqnarray}
When (\ref{5.09})  is valid, and noting  that $\langle \sigma^{\lambda}\rangle^{-(\frac{1}{2}-2\epsilon)}
\langle \sigma_{1}^{\lambda}\rangle^{-(\frac{1}{2}+\epsilon)}\leq
\langle \sigma_{1}^{\lambda}\rangle^{-(\frac{1}{2}-2\epsilon)}\langle \sigma^{\lambda}\rangle^{-(\frac{1}{2}+\epsilon)}$,
we infer that
\begin{eqnarray}
&& K_{3}(\xi_{1},\tau_{1},\xi,\tau)\leq \frac{C|\xi||\xi_{2}|^{-s}}
{\langle\sigma_{1}^{\lambda}\rangle^{\frac{1}{2}-2\epsilon}
\langle\sigma^{\lambda}\rangle^{\frac{1}{2}+\epsilon}\langle\sigma_{2}^{\lambda}
\rangle^{\frac{1}{2}+\epsilon}}\leq
  C\frac{C|\xi|^{2\epsilon}|\xi_{2}|^{\frac{1}{4}-4\epsilon}}
  {\langle\sigma^{\lambda}\rangle^{\frac{1}{2}+\epsilon}\langle\sigma_{2}^{\lambda}
\rangle^{\frac{1}{2}+\epsilon}}.\label{5.013}
\end{eqnarray}
Let $\xi^{\prime}=\xi_{1}$, $\xi_{1}^{\prime}=\xi$,
$\xi_{2}^{\prime}=-\xi_{2}$, $\tau^{\prime}=\tau_{1}$,
 $\tau_{1}^{\prime}=\tau$, $\tau_{2}^{\prime}=-\tau_{2}$,
$L(\xi^{\prime},\tau^{\prime})=F_{1}(\xi_{1},\tau_{1})$,
 $G(\xi_{1}^{\prime},\tau_{1}^{\prime})=F(\xi,\tau)$,
 $H(\xi_{2}^{\prime},\tau_{2}^{\prime})=F_{2}(-\xi_{2},-\tau_{2})$,
$\sigma_{1}^{\prime}=\tau_{1}^{\prime}+\phi^{\lambda}(\xi_{1}^{\prime})=\sigma^{\lambda},
\sigma_{2}^{\prime}=\tau_{2}^{\prime}+\phi^{\lambda}(\xi_{2}^{\prime})=-\sigma_{2}^{\lambda},
\sigma^{\prime}=\tau^{\prime}+\phi^{\lambda}(\xi^{\prime})=\sigma_{1}^{+}.$
Combining (\ref{5.013}) with the above  variable substitution, we have that
\begin{eqnarray}
&& K_{3}(\xi_{1},\tau_{1},\xi,\tau)\leq \frac{C|\xi||\xi_{2}|^{-s}}{\langle\sigma_{1}^{\lambda}
\rangle^{\frac{1}{2}-2\epsilon}\langle\sigma^{\lambda}
\rangle^{\frac{1}{2}+\epsilon}\langle\sigma_{2}^{\lambda}
\rangle^{\frac{1}{2}+\epsilon}}\leq
  C\frac{C|\xi|^{2\epsilon}|\xi_{2}|^{\frac{1}{4}-4\epsilon}}
  {\langle\sigma^{+}\rangle^{\frac{1}{2}+\epsilon}\langle\sigma_{2}^{\lambda}
\rangle^{\frac{1}{2}+\epsilon}}\nonumber\\
&&\leq C
\frac{C|\xi_{1}^{\prime}|^{2\epsilon}|\xi_{2}^{\prime}|^{\frac{1}{4}-4\epsilon}}
  {\langle\sigma_{1}^{\prime}\rangle^{\frac{1}{2}+\epsilon}\langle\sigma_{2}^{\prime}
\rangle^{\frac{1}{2}+\epsilon}}\leq C
 \frac{|\left[\phi^{\lambda}(\xi_{1}^{\prime})\right]^{\prime}-
 \left[\phi^{\lambda}(\xi_{2}^{\prime})\right]^{\prime}|^{\frac{1}{2}}}
 {\langle\sigma_{1}^{\prime}\rangle^{\frac{1}{2}+\epsilon}\langle\sigma_{2}^{\prime}
\rangle^{\frac{1}{2}+\epsilon}} .\label{5.014}
\end{eqnarray}
By using (\ref{5.014}),  the Cauchy-Schwartz inequality,  Plancherel identity and
 the H\"older inequality as well as Lemma 2.9,
we have that
\begin{eqnarray}
&&Int\leq C\int_{\SR^2}\int_{\tiny\begin{array}{cc}\xi^{\prime}
=\sum\limits_{j=1}^{2}\xi_j^{\prime}\\ \tau ^{\prime} =\sum\limits_{j=1}^{2}\tau_j^{\prime}
\end{array}}\frac{|\left[\phi^{\lambda}(\xi_{1}^{\prime})
\right]^{\prime}-\left[\phi^{\lambda}(\xi_{2}^{\prime})\right]^{\prime}|^{\frac{1}{2}}
L(\xi^{\prime},\tau^{\prime})G(\xi_{1}^{\prime},\tau_{1}^{\prime})H(\xi_{2}^{\prime},\tau_{2}^{\prime})}
{\prod\limits_{j=1}^{2}\langle\sigma_{j}^{\prime}\rangle^{\frac{1}{2}+\epsilon}}
d\xi_{1}^{\prime}d\tau_{1}^{\prime}d\xi^{\prime} d\tau^{\prime}\nonumber\\
&&\leq C \|L\|_{L_{\xi\tau}^{2}}
\left\|I^{\frac{1}{2}}\left(\mathscr{F}^{-1}\left(\frac{G}{\langle \sigma_{1}^{\prime}
\rangle ^{\frac{1}{2}+\epsilon}}\right),\mathscr{F}^{-1}\left(\frac{H}{\langle \sigma_{2} ^{\prime}\rangle ^{\frac{1}{2}+\epsilon}}\right)\right)\right\|_{L_{xt}^{2}}\nonumber\\
&&\leq C\|F\|_{L_{\xi\tau}^{2}}\left(\prod_{j=1}^{2}\|F_{j}\|_{L_{\xi\tau}^{2}}\right).\label{5.015}
\end{eqnarray}
When (\ref{5.010})  is valid, this case can be proved similarly to case (\ref{5.09}) of $\Omega_{2}$ of Lemma 5.1.

\noindent (3) Region $\Omega_{3}$. In this subregion, we have that
\begin{eqnarray}
 K_{3}(\xi_{1},\tau_{1},\xi,\tau)\leq \frac{C|\xi||\xi_{1}|^{\frac{3}{2}-4\epsilon}}
 {\langle\sigma^{\lambda}\rangle^{\frac{1}{2}-2\epsilon}\prod\limits_{j=1}^{2}
 \langle\sigma_{j}^{\lambda}\rangle^{\frac{1}{2}+\epsilon}}.\label{5.016}
\end{eqnarray}
Note that $|\sigma^{\lambda}-\sigma_{1}^{\lambda}-\sigma_{2}^{\lambda}|
=\left|3\xi\xi_{1}\xi_{2}-\frac{\xi_{1}^{2}+\xi_{1}\xi_{2}+\xi_{2}^{2}}{\lambda^{4}\xi\xi_{1}\xi_{2}}\right|$.

\noindent Now we consider  (\ref{1.016}), (\ref{1.017}),
separately.

\noindent
When (\ref{1.017}) is valid,
one of (\ref{5.08})-(\ref{5.010}) must occur.

\noindent When (\ref{5.08}) is valid,
we have that
\begin{eqnarray*}
&&K_{3}(\xi_{1},\tau_{1},\xi,\tau)\leq C\frac{|\xi|^{\frac{1}{2}+\epsilon}
|\xi_{1}|^{\frac{1}{2}-4\epsilon}}{\prod\limits_{j=1}^{2}
 \langle\sigma_{j}^{\lambda}\rangle^{\frac{1}{2}+\epsilon}}\nonumber\\
 &&\leq C\frac{|\xi|^{\frac{1}{2}}|\xi_{1}|^{\frac{1}{2}}}{\prod\limits_{j=1}^{2}
 \langle\sigma_{j}^{\lambda}\rangle^{\frac{1}{2}+\epsilon}}\leq C
 \frac{|\xi_{1}^{2}-\xi_{2}^{2}|^{\frac{1}{2}}}{\prod\limits_{j=1}^{2}
 \langle\sigma_{j}^{\lambda}\rangle^{\frac{1}{2}+\epsilon}}\leq
 C\frac{|\left[\phi^{\lambda}(\xi_{1})\right]^{\prime}-
 \left[\phi^{\lambda}(\xi_{2})\right]^{\prime}|^{\frac{1}{2}}}{\prod\limits_{j=1}^{2}
 \langle\sigma_{j}^{\lambda}\rangle^{\frac{1}{2}+\epsilon}}.
\end{eqnarray*}
This case can be proved similarly to (\ref{5.012}).

\noindent When (\ref{5.09}) is valid,
we have that
\begin{eqnarray*}
&&K_{3}(\xi_{1},\tau_{1},\xi,\tau)\leq C\frac{|\xi||\xi_{1}|^{-2s}}
{\langle\sigma^{\lambda}\rangle^{\frac{1}{2}-2\epsilon}\prod\limits_{j=1}^{2}
 \langle\sigma_{j}^{\lambda}\rangle^{\frac{1}{2}+\epsilon}}\leq C
 \frac{|\xi||\xi_{1}|^{-2s}}{\langle\sigma_{1}^{\lambda}\rangle^{\frac{1}{2}-2\epsilon}
 \langle\sigma^{\lambda}\rangle^{\frac{1}{2}+\epsilon}\langle\sigma_{2}^{\lambda}
 \rangle^{\frac{1}{2}+\epsilon}}\nonumber\\&&\leq C\frac{|\xi|^{\frac{1}{2}+2\epsilon}
 |\xi_{1}|^{\frac{1}{2}-8\epsilon}}
 {\langle\sigma^{\lambda}\rangle^{\frac{1}{2}+\epsilon}
 \langle\sigma_{2}^{\lambda}\rangle^{\frac{1}{2}+\epsilon}}\leq C
 \frac{|\xi|^{\frac{1}{2}+2\epsilon}|\xi_{1}|^{\frac{1}{2}-8\epsilon}}
 {\langle\sigma^{\lambda}\rangle^{\frac{1}{2}+\epsilon}
 \langle\sigma_{2}^{\lambda}\rangle^{\frac{1}{2}+\epsilon}}\nonumber\\&&\leq C\frac{|\xi^{2}-\xi_{2}^{2}|^{\frac{1}{2}}}{\langle\sigma^{\lambda}\rangle^{\frac{1}{2}+\epsilon}
 \langle\sigma_{2}^{\lambda}\rangle^{\frac{1}{2}+\epsilon}}\leq C\frac{|\left[\phi^{\lambda}(\xi)\right]^{\prime}-\left[\phi^{\lambda}(\xi_{2})\right]^{\prime}|^{\frac{1}{2}}
 }{\langle\sigma^{\lambda}\rangle^{\frac{1}{2}+\epsilon}
 \langle\sigma_{2}^{\lambda}\rangle^{\frac{1}{2}+\epsilon}}.
\end{eqnarray*}
This case can be proved similarly to  (\ref{5.015}).

\noindent When (\ref{5.010})  is valid,
this case can be proved similarly to case (\ref{5.09}) of $\Omega_{3}$ of Lemma 5.1.

\noindent When (\ref{1.016}) is valid,
we have   $|\xi|\sim |\xi_{1}|^{-1}$. Thus
\begin{eqnarray}
 &&K_{3}(\xi_{1},\tau_{1},\xi,\tau)\leq \frac{C|\xi||\xi_{1}|^{\frac{3}{2}-12\epsilon}}
 {\langle\sigma^{\lambda}\rangle^{\frac{1}{2}-2\epsilon}
 \prod\limits_{j=1}^{2}\langle\sigma_{j}^{\lambda}\rangle^{\frac{1}{2}+\epsilon}}
 \leq \frac{C|\xi_{1}|^{\frac{1}{2}-12\epsilon}}
 {\langle\sigma^{\lambda}\rangle^{\frac{1}{2}-2\epsilon}\prod\limits_{j=1}^{2}
 \langle\sigma_{j}^{\lambda}\rangle^{\frac{1}{2}+\epsilon}}
 .\label{5.017}
\end{eqnarray}
Let $\xi^{\prime}=\xi_{1}$, $\xi_{1}^{\prime}=\xi$,
$\xi_{2}^{\prime}=-\xi_{2}$, $\tau^{\prime}=\tau_{1}$,
 $\tau_{1}^{\prime}=\tau$, $\tau_{2}^{\prime}=-\tau_{2}$,
$L(\xi^{\prime},\tau^{\prime})=F_{1}(\xi_{1},\tau_{1})$,
 $G(\xi_{1}^{\prime},\tau_{1}^{\prime})=F(\xi,\tau)$,
 $H(\xi_{2}^{\prime},\tau_{2}^{\prime})=F_{2}(-\xi_{2},-\tau_{2})$,
$\sigma_{1}^{\prime}=\tau_{1}^{\prime}+\phi^{\lambda}(\xi_{1}^{\prime})=\sigma^{\lambda},
\sigma_{2}^{\prime}=\tau_{2}^{\prime}+\phi^{\lambda}(\xi_{2}^{\prime})=-\sigma_{2}^{\lambda},
\sigma^{\prime}=\tau^{\prime}+\phi(\xi^{\prime})=\sigma_{1}^{\lambda}.$
Combining (\ref{5.017}) with the above  variable substitution, we have that
\begin{eqnarray}
&& K_{3}(\xi_{1},\tau_{1},\xi,\tau)\leq \frac{C|\xi_{1}|^{\frac{1}{2}-12\epsilon}}
{\langle\sigma^{\lambda}\rangle^{\frac{1}{2}-2\epsilon}
\prod\limits_{j=1}^{2}\langle\sigma_{j}^{\lambda}\rangle^{\frac{1}{2}+\epsilon}}\leq
 \frac{C|\xi^{\prime}|^{\frac{1}{2}-12\epsilon}}{\langle\sigma_{1}^{\prime}
\rangle^{\frac{1}{2}-2\epsilon}\langle\sigma^{\prime}\rangle^{\frac{1}{2}+\epsilon}
\langle\sigma_{2}^{\prime}\rangle^{\frac{1}{2}+\epsilon}} \nonumber\\
&&\leq C\frac{|\left[\phi^{\lambda}(\xi_{1}^{\prime})
\right]^{\prime}-\left[\phi^{\lambda}(\xi_{2}^{\prime})\right]^{\prime}|^{\frac{1}{4}-6\epsilon}}
{\langle\sigma^{\prime}\rangle^{\frac{1}{2}+\epsilon}\langle\sigma_{1}^{\prime}\rangle^{\frac{1}{2}-2\epsilon}
\langle\sigma_{2}^{\prime}\rangle^{\frac{1}{2}+\epsilon}}.\label{5.018}
\end{eqnarray}
By using (\ref{5.018}), together with  the Cauchy-Schwartz inequality,  Plancherel identity and
 the H\"older inequality as well as Lemma 2.9,
we infer that
\begin{eqnarray}
&&\hspace{-0.5cm}Int\leq C\int_{\SR^2}\int_{\tiny\begin{array}{cc}\xi^{\prime}
=\sum\limits_{j=1}^{2}\xi_j^{\prime}\\ \tau ^{\prime} =\sum\limits_{j=1}^{2}
\tau_j^{\prime}\end{array}}\frac{|\left[\phi^{\lambda}(\xi_{1}^{\prime})
\right]^{\prime}-\left[\phi^{\lambda}(\xi_{2}^{\prime})\right]^{\prime}|^{\frac{1}{4}-6\epsilon}
L(\xi^{\prime},\tau^{\prime})G(\xi_{1}^{\prime},\tau_{1}^{\prime})H(\xi_{2}^{\prime},\tau_{2}^{\prime})}
{\langle\sigma^{\prime}\rangle^{\frac{1}{2}+\epsilon}\langle\sigma_{1}^{\prime}\rangle^{\frac{1}{2}-2\epsilon}
\langle\sigma_{2}^{\prime}\rangle^{\frac{1}{2}+\epsilon}}
d\xi_{1}^{\prime}d\tau_{1}^{\prime}d\xi^{\prime} d\tau^{\prime}\nonumber\\
&&\leq C\left \|\frac{L}{\langle\sigma^{\prime}\rangle^{\frac{1}{2}+\epsilon}}\right\|_{L_{\xi\tau}^{2}}
\left\|I^{\frac{1}{2}}\left(\mathscr{F}^{-1}\left(\frac{G}{\langle \sigma_{1}^{\prime}
\rangle ^{\frac{1}{2}-2\epsilon}}\right),\mathscr{F}^{-1}\left(\frac{H}{\langle \sigma_{2} ^{\prime}\rangle ^{\frac{1}{2}+\epsilon}}\right)\right)\right\|_{L_{xt}^{2}}\nonumber\\
&&\leq C\|F\|_{L_{\xi\tau}^{2}}\left(\prod_{j=1}^{2}\|F_{j}\|_{L_{\xi\tau}^{2}}\right).\label{5.019}
\end{eqnarray}

\noindent (4) Region $\Omega_{4}$. In this subregion, we have that
\begin{eqnarray}
 K_{3}(\xi_{1},\tau_{1},\xi,\tau)\leq \frac{C|\xi_{1}|^{\frac{7}{4}-6\epsilon}}
 {\langle\sigma^{\lambda}\rangle^{\frac{1}{2}-2\epsilon}
 \prod\limits_{j=1}^{2}\langle\sigma_{j}^{\lambda}\rangle^{\frac{1}{2}+\epsilon}}.\label{5.020}
\end{eqnarray}
In this case, (\ref{1.017}) is valid. Thus, one of  (\ref{5.08})-(\ref{5.010}) must occur.

\noindent
When (\ref{5.08}) is valid, we get from (\ref{5.020}) that
\begin{eqnarray}
 K_{3}(\xi_{1},\tau_{1},\xi,\tau)\leq \frac{C|\xi_{1}|^{\frac{1}{4}}}{
 \prod\limits_{j=1}^{2}\langle\sigma_{j}^{\lambda}\rangle^{\frac{1}{2}+\epsilon}}
 \leq C\frac{\prod\limits_{j=1}^{2}|\xi_{j}|^{\frac{1}{8}}}{
 \prod\limits_{j=1}^{2}\langle\sigma_{j}^{\lambda}\rangle^{\frac{1}{2}+\epsilon}}.
 \label{5.021}
\end{eqnarray}
By using (\ref{5.021}), the Cauchy-Schwartz inequality, the Plancherel identity and
 the H\"older inequality as well as (\ref{2.022}),
we conclude  that
\begin{eqnarray}
&&Int\leq C\int_{\SR^2}\int_{\tiny\begin{array}{cc}\xi=\sum\limits_{j=1}^{2}\xi_j\\
\tau =\sum\limits_{j=1}^{2}\tau_j\end{array}}\frac{
F(\xi,\tau)\prod\limits_{j=1}^{2}|\xi_{j}|^{\frac{1}{8}}F_{j}(\xi_{j},\tau_{j})}{
\prod\limits_{j=1}^{2}\langle\sigma_{j}^{\lambda}\rangle^{\frac{1}{2}+\epsilon}}
d\xi_{1}d\tau_{1}d\xi d\tau\nonumber\\
&&\leq C \left\|F\right\|_{L_{\xi\tau}^{2}}\left(\prod_{j=1}^{2}
\left\|D_{x}^{\frac{1}{8}}P^{4}\mathscr{F}^{-1}\left(\frac{F_{j}}{\langle \sigma_{j}^{\lambda} \rangle ^{\frac{1}{2}+\epsilon}}\right)\right\|_{L_{xt}^{4}}\right)\nonumber\\
&&\leq C\|F\|_{L_{\xi\tau}^{2}}\left(\prod_{j=1}^{2}\|F_{j}\|_{L_{\xi\tau}^{2}}\right).\label{5.022}
\end{eqnarray}
When (\ref{5.09}) is valid, we infer  from (\ref{5.020}) that
\begin{eqnarray}
 K_{3}(\xi_{1},\tau_{1},\xi,\tau)\leq \frac{C|\xi_{1}|^{\frac{1}{4}}}
 {\langle\sigma^{\lambda}\rangle^{\frac{1}{2}-2\epsilon}
 \langle\sigma_{2}^{\lambda}\rangle^{\frac{1}{2}+\epsilon}}\leq\
  C\frac{\prod\limits_{j=1}^{2}|\xi_{j}|^{\frac{1}{8}}}
 {\langle\sigma^{\lambda}\rangle^{\frac{1}{2}-2\epsilon}
 \langle\sigma_{2}^{\lambda}\rangle^{\frac{1}{2}+\epsilon}}.\label{5.023}
\end{eqnarray}
Let $\xi^{\prime}=\xi_{1}$, $\xi_{1}^{\prime}=\xi$,
$\xi_{2}^{\prime}=-\xi_{2}$, $\tau^{\prime}=\tau_{1}$,
 $\tau_{1}^{\prime}=\tau$, $\tau_{2}^{\prime}=-\tau_{2}$,
$L(\xi^{\prime},\tau^{\prime})=F_{1}(\xi_{1},\tau_{1})$,
 $G(\xi_{1}^{\prime},\tau_{1}^{\prime})=F(\xi,\tau)$,
 $H(\xi_{2}^{\prime},\tau_{2}^{\prime})=F_{2}(-\xi_{2},-\tau_{2})$,
$\sigma_{1}^{\prime}=\tau_{1}^{\prime}+\phi^{\lambda}(\xi_{1}^{\prime})=\sigma^{\lambda},
\sigma_{2}^{\prime}=\tau_{2}^{\prime}+\phi^{\lambda}(\xi_{2}^{\prime})=-\sigma_{2}^{\lambda},
\sigma^{\prime}=\tau^{\prime}+\phi(\xi^{\prime})=\sigma_{1}^{\lambda}.$
Combining (\ref{5.023}) with the above  variable substitution, we have that
\begin{eqnarray}
&& K_{3}(\xi_{1},\tau_{1},\xi,\tau)\leq \frac{C\prod\limits_{j=1}^{2}
|\xi_{j}^{\prime}|^{\frac{1}{8}}}
{\langle\sigma_{1}^{\prime}\rangle^{\frac{1}{2}-2\epsilon}
\langle\sigma_{2}^{\prime}\rangle^{\frac{1}{2}+\epsilon}
} .\label{5.024}
\end{eqnarray}
Combining (\ref{5.024}), (\ref{2.022}) with a proof similar to (\ref{5.022}),
 and noting that
$\frac{3}{4}\left(\frac{1}{2}+\epsilon\right)<\frac{1}{2}-2\epsilon,$
 we obtain  that
\begin{eqnarray}
Int\leq C\|F\|_{L_{\xi\tau}^{2}}\left(\prod_{j=1}^{2}
\|F_{j}\|_{L_{\xi\tau}^{2}}\right)
.\label{5.025}
\end{eqnarray}
\noindent When (\ref{5.010}) is valid, this case can be proved similarly to
 (\ref{5.09}) of $\Omega_{4}$ in Lemma 5.1.

\noindent (5) Region $\Omega_{5}$. This case
can be proved similarly to Subregion (4) of  Lemma 5.1.

This completes the proof of Lemma 5.1, and thus also the
proof of Theorem 1.2.

\noindent {\bf Remark 8.}As we mentioned earlier,  the bilinear estimate in Theorem 1.2
 leads to the proof of  the local
well-posedness of the Ostrovsky equation  in $H^{s}(\R)$ for
 $ s>-\frac{3}{4}$ with the
 aid of Lemma 2.8 and a fixed point argument.

\bigskip
\bigskip

\leftline{\large \bf Acknowledgments}
\noindent

 This work was supported by the NSF grant 1620449 and the
 Natural Science Foundation of China
 under grant  11401180.  The first author is also supported
 by the Young Core Teachers
  Program of Henan Normal University and Henan Province Research Fund 15A110033.

\end{document}